\RequirePackage{ifthen}
\newboolean{MOR}
\setboolean{MOR}{false}

\ifthenelse {\boolean{MOR}}
{
\documentclass[moor]{informs1}
} {
\documentclass[10pt]{article}
\usepackage[margin=1in]{geometry}
}

\usepackage{graphicx}
\usepackage{subfigure,epsfig,url,psfrag}
\usepackage{cancel}
\usepackage{amsfonts}
\usepackage{amsmath}
\usepackage{amssymb}
\usepackage{verbatim,booktabs}
\usepackage{float}
\usepackage{mathrsfs}
\usepackage{colonequals}
\usepackage{color}

\usepackage{amsopn}

\ifthenelse {\boolean{MOR}}
{
\usepackage{natbib}
 \NatBibNumeric
 \bibpunct[, ]{[}{]}{,}{n}{}{,}%

\usepackage[colorlinks=true,breaklinks=true,bookmarks=true,urlcolor=blue,
     citecolor=blue,linkcolor=blue,bookmarksopen=false,draft=false]{hyperref}

\def\EMAIL#1{\href{mailto:#1}{#1}}

\TheoremsNumberedThrough     
\EquationsNumberedThrough    

\newcommand{\sctn}[1]{\section{#1.}}
\newcommand{\sbsctn}[1]{\subsection{#1.}}
\newcommand{\sbsbsctn}[1]{\subsubsection{#1.}}

\newenvironment{prf}[1][]
{\begin{proof}{Proof.}}
{\hfill\ensuremath{\square} \end{proof} \smallskip}
}
{
\usepackage{hyperref}
\usepackage{amsthm}

\newtheorem{theorem}{Theorem}
\newtheorem{lemma}{Lemma}
\newtheorem{proposition}{Proposition}

\newcommand{\sctn}[1]{\section{#1}}
\newcommand{\sbsctn}[1]{\subsection{#1}}
\newcommand{\sbsbsctn}[1]{\subsubsection{#1}}

\newenvironment{prf}[1][]
{\begin{proof}}
{\end{proof}}
}

\DeclareMathOperator{\prob}{\mathbb{P}}
\DeclareMathOperator{\avg}{\mathbb{E}}

\DeclareMathOperator{\reg}{reg}
\DeclareMathOperator{\sub}{s.t.}
\DeclareMathOperator{\find}{find}

\newcommand{\din}{d_{\text{in}}}
\newcommand{\dout}{d_{\text{out}}}
\newcommand{\Z}{\mathbb Z}
\newcommand{\R}{\mathbb R}

\newcommand{\tin}{\mathrm{in}}
\newcommand{\tout}{\mathrm{out}}
\newcommand{\con}{\mathrm{con}}
\renewcommand{\deg}{\mathrm{deg}}

\newcommand{\ie}{i.e., }

\newcommand{\G}{\mathcal G}
\newcommand{\I}{\mathcal I}

\newcommand{\tavg}{\mathrm{avg}}

\begin{document}

\newcommand{\fndng}{A.~Del~Pia is partially funded by  ONR grant N00014-19-1-2322. Any opinions, findings, and conclusions or recommendations expressed in this material are those of the authors and do not necessarily reflect the views of the Office of Naval Research. D.~Kunisky is supported by ONR Award N00014-20-1-2335, a Simons Investigator Award to Daniel Spielman, and NSF grants DMS-1712730 and DMS-1719545. Part of this work was performed while D.~Kunisky was with New York University.}

\newcommand{\ttl}{Linear Programming and Community Detection\thanks{\fndng}}

\newcommand{\bstrct}{The problem of community detection with two equal-sized communities is closely related to the minimum graph bisection problem over certain random graph models. In the stochastic block model distribution over networks with community structure, a well-known semidefinite programming (SDP) relaxation of the minimum bisection problem recovers the underlying communities whenever possible. Motivated by their superior scalability,
we study the theoretical performance of linear programming (LP) relaxations of the minimum bisection problem for the same random models.
We show that, unlike the SDP relaxation that undergoes a phase transition in the logarithmic average degree regime, the LP relaxation fails in recovering the planted bisection with high probability in this regime.  We show that the LP relaxation instead exhibits a transition from recovery to non-recovery in the linear average degree regime. Finally, we present non-recovery conditions for graphs with average degree strictly between linear and logarithmic.}


\newcommand{\kywrds}{Community detection, minimum bisection problem, linear programming, metric polytope}

\ifthenelse {\boolean{MOR}}
{
\TITLE{\ttl}
\RUNTITLE{\ttl}

\ARTICLEAUTHORS{
\AUTHOR{Alberto Del Pia}
\AFF{Department of Industrial and Systems Engineering \& Wisconsin Institute for Discovery, University of Wisconsin-Madison, \EMAIL{delpia@wisc.edu}}
\AUTHOR{Aida Khajavirad}
\AFF{Department of Industrial \& Systems Engineering, Lehigh University, \EMAIL{aida@lehigh.edu}}
\AUTHOR{Dmitriy Kunisky}
\AFF{Department of Computer Science, Yale University, \EMAIL{dmitriy.kunisky@yale.edu}}}
\RUNAUTHOR{Del Pia et al.}

\ABSTRACT{\bstrct}
\KEYWORDS{\kywrds}
}
{
\title{\ttl}

\author{
Alberto Del Pia
\thanks{Department of Industrial and Systems Engineering \& Wisconsin Institute for Discovery, University of Wisconsin-Madison.
E-mail: {\tt delpia@wisc.edu}}
\and
Aida Khajavirad
\thanks{Department of Industrial and System Engineering, Lehigh University.
  E-mail: {\tt aida@lehigh.edu}.}
\and
Dmitriy Kunisky
\thanks{Department of Computer Science, Yale University.
E-mail: {\tt  dmitriy.kunisky@yale.edu}}
}
}

\date{May 11, 2022}

\maketitle



\ifthenelse {\boolean{MOR}}
{}
{
\begin{abstract}
\bstrct
\end{abstract}

\emph{Key words:}
\kywrds
}

\sctn{Introduction}
\label{sec: intro}

Performing \emph{community detection} or \emph{graph clustering} in large networks is a central problem in applied disciplines including biology, social sciences, and engineering.
In community detection, we are given a network of nodes and edges, which may represent anything from social actors and their interactions, to genes and their functional cooperation, to circuit components and their physical connections.
We then wish to find \emph{communities}, or subsets of nodes that are densely connected to one another.
The exact solution of such problems typically amounts to solving NP-hard graph partitioning problems; hence, practical techniques instead produce approximations based on various heuristics~\cite{SnijNow97,CarImp01,BicChen09,Fort10}.

\paragraph{The stochastic block model.} Various \emph{generative models} of random networks with community structure have been proposed as a simple testing ground for the numerous available algorithmic techniques.
Analyzing the performance of algorithms in this way has the advantage of not relying on individual test cases from particular domains, and of capturing the performance on \emph{typical} random problem instances, rather than worst-case instances where effective guarantees of performance can seldom be made.
The stochastic block model (SBM) is the most widely-studied generative model for community detection.
Under the SBM, nodes are assigned to one of several communities (the ``planted'' partition or community assignment), and two nodes are connected with a probability depending only on their communities.
In the simplest case, if there is an even number $n$ of nodes, then we may assign the nodes to two communities of size $n / 2$, where nodes in the same community are connected with probability $p = p(n)$ and nodes in different communities are connected with probability $q = q(n)$, for some $p > q$.
This is the so-called \emph{symmetric assortative} SBM with two communities.

In this paper, we study the problem of \emph{exact recovery} (henceforth simply \emph{recovery}) of the communities under this model.
That is, we are interested in algorithms that, with high probability (\ie probability tending to 1 as $n \to \infty$), recover the assignment of nodes to communities correctly.
For this task, the estimator that is most likely to succeed is the maximum \emph{a posteriori} (MAP) estimator, which coincides with the \emph{minimum bisection} of the graph, the assignment with the least number of edges across communities (see Chapter~3 of~\cite{Abbe18}).
Computing the minimum bisection is NP-hard \cite{GarJoSto74}, and the best known polynomial-time approximations have a poly-logarithmic worst-case multiplicative error in the size of the bisection~\cite{KraFei06}.

\paragraph{Semidefinite programming relaxations.} However, there is still hope to solve the recovery problem under the SBM efficiently, since that only requires an algorithm to perform well on typical random graphs.
Indeed, in this setting, one recent stream of research has shown that semidefinite programming (SDP) relaxations of various community detection and graph clustering problems successfully recover communities under suitable generative models~\cite{AbbBanHal16,MixVilWar17,IduMixPetVil17,Ban18,LinStr19,LiLiLi20}.
Generally speaking, these works first provide deterministic sufficient conditions for a given community assignment to be the unique optimal solution of the SDP relaxation, and then show that those conditions hold with high probability under a given model.
For the symmetric assortative SBM with two communities, the work \cite{AbbBanHal16} gave a tight characterization of the values of $p$ and $q$ for which \emph{any} algorithm, regardless of computational complexity, can recover the community assignments (often called an ``information-theoretic'' threshold).
In particular, the authors found that recovery changes from possible to impossible in the asymptotic regime $p(n), q(n) \sim \log n / n$.
The authors also conjectured that an SDP relaxation in fact achieves this limit, and proved a partial result in this direction, later extended to a full proof of the conjecture by \cite{HajWuXu16,Ban18}.
These results indicate that, remarkably, the polynomial-time SDP relaxation succeeds in recovery whenever the (generally intractable) MAP estimator does, in a suitable asymptotic sense.

\paragraph{Linear programming relaxations.} It is widely accepted that for problems of comparable size, state-of-the-art solvers for linear programming (LP) significantly outperform those for SDP in both speed and scalability.
Yet, in contrast to the rich literature on SDP relaxations, LP relaxations for community detection have not been studied.
In this paper, motivated by their desirable practical properties, we study the theoretical performance of LP relaxations of the minimum bisection problem for recovery under the SBM. Recently, in~\cite{AntoAida20,AidaTonio21,dPMa21m,dPIda21}, the authors investigate the recovery properties of LP relaxations for K-means clustering, join object matching, K-medians clustering, and Boolean tensor factorization, respectively.

\paragraph{Integrality gaps of LP relaxations for graph cut problems.}
Perhaps the most similar line of prior work concerns LP relaxations for the maximum cut and the sparsest cut problems.
Poljak and Tuza~\cite{PolTuz94} consider a well-known LP relaxation of maximum cut problem often referred to as the metric relaxation~\cite{dl:97}.
They show that for Erd\H{o}s-R\'{e}nyi graphs with edge probability $O(\mathsf{polylog}(n) / n)$, the metric relaxation yields a trivial integrality gap of $2$, while for denser graphs with edge probability $\Omega(\sqrt{\log n / n})$, this LP provides smaller integrality gaps.
In~\cite{AviUme03}, the authors consider an LP relaxation obtained by adding a large class of inequalities to the metric relaxation.
Yet, this stronger LP does not improve the trivial integrality gap of the metric relaxation for edge probability $O(\mathsf{polylog}(n) / n)$.
In~\cite{VegKen07}, the authors show that, for high-girth graphs, for any fixed $k$, the LP relaxation obtained after $k$ rounds of the Sherali-Adams hierarchy does not improve the trivial integrality gap either. In fact, the authors of~\cite{ChaMakMak09} prove that, for random $d$-regular graphs, for every $\epsilon > 0$, there exists $\gamma = \gamma(d, \epsilon) > 0$ such that integrality gap of the LP relaxation obtained after $n^{\gamma}$ rounds of Sherali-Adams is $2 - \epsilon$.
In contrast to all of these negative results, the recent works~\cite{DonSch18,HopSchTre19} provide a partial ``redemption'' of the Sherali-Adams relaxation for approximating the maximum cut, by showing that LP relaxations obtained from $n^{O(1)}$ rounds of Sherali-Adams obtain non-trivial worst-case approximation ratios.
For the sparsest cut problem, the authors of \cite{LeiRao99} give an $O(\log n)$-approximation algorithm based on the metric relaxation, and show that constant degree expander graphs in fact yield a matching integrality gap of $\Theta(\log n)$ for this relaxation.

We draw attention to the unifying feature that, in all of these results, the upper bounds produced by LPs for dense random graphs are significantly better than those for sparser random graphs.
  Our results describe another setting, concerning performance on the statistical task of recovering a planted bisection rather than approximating the size of the minimum bisection, where the same phenomenon occurs.

\paragraph{Other algorithmic approaches.} Recent results have also shown that spectral methods can achieve optimal performance in the SBM as well~\cite{Abbe18,AbbFanWanZho20}. As these methods only require estimating the leading eigenvector of a matrix, they are faster than both LP and SDP relaxations. However, it is still valuable to study the theoretical properties of LP relaxations for community detection. First, LP methods provide a general strategy for approximating a wide range of combinatorial graph problems, and understanding the behavior of LP relaxations for the minimum bisection problem sheds light on their applicability to other situations where similar spectral methods do not exist. Second, works like \cite{MoiPerWei16,RicJavMon16} suggest that convex relaxations enjoy robustness to adversarial corruptions of the inputs for statistical problems that spectral methods do not, making it of practical interest to understand the most efficient convex relaxation algorithms for solving statistical problems like recovery in the SBM.

\paragraph{Our contribution.} Our work serves as the first average-case analysis of recovery properties of LP relaxations for the minimum bisection problem under the SBM.
Proceeding similarly to the works on SDP relaxations described above, we begin with a deterministic analysis.
First, we obtain necessary and sufficient conditions for a planted bisection in a graph to be the unique minimum bisection
(see Theorems~\ref{th: IP} and~\ref{th: IP tight}).
These conditions are given in terms of certain simple measures of within-cluster and inter-cluster connectivity.
Next, we consider an LP relaxation of
the minimum bisection problem obtained by outer-approximating the cut-polytope by the
metric-polytope~\cite{dl:97}. Under certain regularity assumptions on the input graph, we derive a
sufficient condition under which the LP recovers the planted bisection (see Propositions~\ref{optimal} and~\ref{uniqueness}).
This condition is indeed tight in the worst-case.
Moreover, we give an extension of this result to general (irregular) graphs, provided that regularity can be achieved by adding and removing edges in an appropriate way (see Theorem~\ref{deter}).
Finally, we present a necessary condition for recovery using the LP relaxation~(see Theorem~\ref{lp-nonrecovery-cond}).
This condition depends on how the average distance between pairs of nodes in the graph compares to the maximum such distance, the \emph{diameter} of the graph. These two parameters have favorable properties for our subsequent probabilistic analysis.

Next, utilizing our deterministic sufficient and necessary conditions, we perform a probabilistic analysis under the SBM.
Namely, we show that if $p = p(n)$ and $q = q(n)$ are constants independent of $n$, then the LP recovers the planted bisection with high probability, provided that $q \leq p - \frac{1}{2}$~(see Theorem~\ref{SBMdense}).
Conversely, if $q > \max\{2p-1, \frac{1}{2} (3-p-\sqrt{(3-p)^2-4p})\}$, then the LP fails to recover the planted bisection
with high probability~(see Theorem~\ref{noRecoveryDense}.1).
Thus, within the \emph{very dense} asymptotic regime $p(n), q(n) = \Theta(1)$, there are some parameters for which the LP achieves recovery, and others for which it does not.
On the other hand, in the \emph{logarithmic} regime $p(n), q(n) = \Theta(\log n / n)$, we prove that with high probability the LP fails to recover the planted bisection~(see Theorem~\ref{noRecoverySparse}).
We also present a collection of non-recovery conditions for the SBM in between these two regimes; that is, when
$p(n), q(n) = \Theta(n^{-\omega})$ for some $0 <\omega < 1$~(see Theorem~\ref{noRecoveryDense}.2-3).
In summary, we find that LP relaxations do not have the desirable theoretical properties of their SDP counterparts under the SBM, but rather have a novel transition between recovery and non-recovery in a different asymptotic regime.

\paragraph{Outline.} The remainder of the paper is organized as follows.
In Section~\ref{sec: bisection} we consider the minimum bisection problem
and obtain necessary and sufficient conditions for recovery.
Subsequently, we consider the LP relaxation in Section~\ref{sec: lp} and obtain necessary conditions and sufficient conditions under which the planted bisection is the unique optimal solution of this LP. In Section~\ref{sec: SBM} we address the question of recovery under the SBM in various regimes.
Some technical results that are omitted in the previous sections are provided in Section~\ref{sec:proofs}.

\sctn{The minimum bisection problem}
\label{sec: bisection}

Let $G = (V,E)$ be a graph.
A \emph{bisection of $V$} is a partition of $V$ into two subsets of equal cardinality.
Clearly a bisection of $V$ only exists if $|V|$ is even.
The \emph{cost} of the bisection is the number of edges connecting the two sets.
The \emph{minimum bisection problem} is the problem of finding a bisection of minimum cost in a given graph.
This problem is known to be NP-hard~\cite{GarJoSto74} and it is a prototypical graph partitioning problem which arises in numerous applications.
Throughout the paper, we assume $V = \{1,\dots, n\}$.
Furthermore, we denote by $(i,j)$ an edge in $E$ with ends $i<j$.

Let $G$ be a given graph and assume that there is also a fixed bisection of $V$ which is unknown to us.
We refer to this fixed bisection as the \emph{planted bisection}, and we denote it by $V_1,V_2$.
The question that we consider in this section is the following: When does solving the minimum bisection problem recover the planted bisection? Clearly this can only happen when the planted bisection is the unique optimal solution.
From now on, we say that an optimization problem \emph{recovers the planted bisection} if its unique optimal solution corresponds to the planted bisection.



We present a necessary and sufficient condition under which the minimum bisection problem recovers the planted bisection.
This condition is expressed in terms of two parameters $\din,\dout$ of the given graph $G$.
Roughly speaking, $\din$ is a measure of \emph{intra-cluster connectivity} while $\dout$ is a measure of \emph{inter-cluster connectivity}.
Such parameters are commonly used in order to obtain performance guarantees for various clustering-type algorithms
(see for example~\cite{KumKan10,LiLiLi20,LinStr19}).
These quantities will also appear naturally in the course of our construction of a dual certificate for the LP relaxation in Theorem~\ref{deter}.
In order to formally define these two parameters, we recall that the \emph{subgraph} of $G$ induced by a subset $U$ of $V$, denoted by $G[U]$, is the graph with node set $U$, and its edge set is the subset of all the edges in $E$ with both ends in $U$.
Furthermore, we denote by $E_1$ and $E_2$ the edge sets of $G[V_1]$ and $G[V_2]$, respectively, and we
define $G_0$ as the graph $G_0 = (V, E_0)$, where $E_0 := E \setminus (E_1 \cup E_2)$.
The parameter $\din$ is then defined as the minimum degree of the nodes of $G[V_1] \cup G[V_2]$, while the parameter $\dout$ is the maximum degree of the nodes of $G_0$.

\sbsctn{A sufficient condition for recovery}

We first present a condition that guarantees that the minimum bisection problem recovers the planted bisection.

\begin{theorem}
\label{th: IP}
If $\din - \dout > n/4 - 1$, then the minimum bisection problem recovers the planted bisection.
\end{theorem}

\begin{prf}
We prove the contrapositive, thus we assume that the minimum bisection problem does not recover the planted bisection, and we show $\din - \dout \le n/4 - 1$.
Let $\tilde V_1, \tilde V_2$ denote the optimal solution to the minimum bisection problem.

Let $U_1 := V_1 \cap \tilde V_1$.
Up to switching $\tilde V_1$ and $\tilde V_2$, we can assume without loss of generality that $|U_1| \le n/4$.
Since the minimum bisection problem does not recover the planted bisection, we have that $U_1$ is nonempty, and we denote by $t$ its cardinality, i.e., $t := |U_1|$.
Let $U_2 := V_2 \cap \tilde V_2$ and note that $|U_2| = t$.
We also define $W_1 := V_1 \cap \tilde V_2$ and $W_2 := V_2 \cap \tilde V_1$.
In particular, we have $V_1 = U_1 \cup W_1$, $V_2 = U_2 \cup W_2$, $\tilde V_1 = U_1 \cup W_2$, and $\tilde V_2 = U_2 \cup W_1$.
See Fig.~\ref{fig: flip} for an illustration.
\begin{figure}[htb]
\begin{center}
\includegraphics[width=0.9\textwidth]{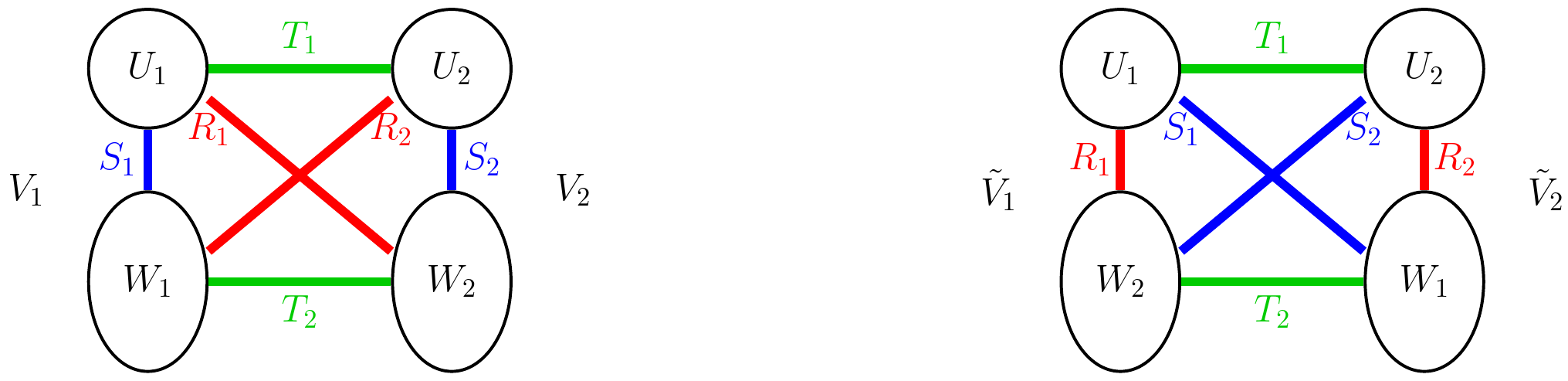}
\end{center}
\caption{Both graphs in the figure are $G$.
The graph on the left represents the planted bisection $V_1,V_2$.
The graph on the right represents the optimal solution of the minimum bisection problem $\tilde V_1,\tilde V_2$.}
\label{fig: flip}
\end{figure}

We define the following sets of edges:
\begin{align*}
R_1 := & \{ e \in E : \text{$e$ has one endnode in $U_1$ and one endnode in $W_2$}\}, \\
R_2 := & \{ e \in E : \text{$e$ has one endnode in $W_1$ and one endnode in $U_2$}\}, \\
S_1 := & \{ e \in E : \text{$e$ has one endnode in $U_1$ and one endnode in $W_1$}\}, \\
S_2 := & \{ e \in E : \text{$e$ has one endnode in $U_2$ and one endnode in $W_2$}\}, \\
T_1 := & \{ e \in E : \text{$e$ has one endnode in $U_1$ and one endnode in $U_2$}\}, \\
T_2 := & \{ e \in E : \text{$e$ has one endnode in $W_1$ and one endnode in $W_2$}\}.
\end{align*}
Furthermore, let $R := R_1 \cup R_2$, $S := S_1 \cup S_2$, and $T := T_1 \cup T_2$.

The optimality of the bisection $\tilde V_1, \tilde V_2$ implies that $|S| + \cancel{|T|} \le |R| + \cancel{|T|}$, thus $|S| \le |R|$.
Next, we obtain an upper bound on $|R|$ and a lower bound on $|S|$.
By definition of $\dout$, we have $|R_1| \le t \dout$ and $|R_2| \le t \dout$, thus $|R| \le 2t \dout$.
Consider now the set $S_1$.
The sum of the degrees of the nodes in $U_1$ in the graph $G[V_1]$ is at least $t \din$, while the sum of the degrees of the nodes of the graph $G[U_1]$ is at most $t(t-1)$.
Thus we have $|S_1| \ge t(\din - t + 1)$.
Symmetrically, $|S_2| \ge t(\din - t + 1)$, thus $|S| \ge 2 t(\din - t + 1)$.
We obtain
\begin{align*}
2t (\din - t + 1) \le |S| \le |R| \le 2t \dout
\quad \Rightarrow \quad
\din - t + 1 \le \dout.
\end{align*}
Since $t \le n/4$, we derive
$\din - \dout \le n/4 - 1$.
\end{prf}

\sbsctn{A necessary condition for recovery}

Next, we present a necessary condition for recovery of the planted bisection.
To this end, we make use of the following graph theoretic lemma.


\begin{lemma}
\label{lem: regular}
For every positive even $t$, and every $d$ with $0 \le d \le t/2$, there exists a $d$-regular bipartite graph with $t$ nodes and where each set in the bipartition contains $t/2$ nodes.
Furthermore, for every positive even $t$, and every $d$ with $0 \le d \le t-1$, there exists a $d$-regular graph with $t$ nodes.
\end{lemma}

\begin{prf}
Fix any positive even $t$, and let $U_1, U_2$ be two disjoint sets of nodes of cardinality $t/2$.
We start by proving the first part of the statement.
For every $d$ with $0 \le d \le t/2$ we explain how to construct a $d$-regular bipartite graph $B_d$ with bipartition $U_1, U_2$.
The graph $B_{t/2}$ is the complete bipartite graph.
We recursively define $B_{d-1}$ from $B_d$ for every $1 \le d \le t/2$.
The graph $B_d$ is regular bipartite of positive degree and so it has a perfect matching (see Corollary 16.2b in \cite{SchBookCO}), which we denote by $M_d$.
The graph $B_{d-1}$ is then obtained from $B_d$ by removing all the edges in $M_d$.

Next we prove the second part of the statement.
Due to the first part of the proof, we only need to construct a $d$-regular graph $G_d$ with node set $U_1 \cup U_2$, for every $d$ with $t/2 + 1 \le d \le t-1$.
Let $H$ be the graph with node set $U_1 \cup U_2$ and whose edge set consists of all the edges with both ends in the same $U_i$, $i=1,2$.
Note that $H$ is a $t/2-1$ regular graph.
For every $t/2 + 1 \le d \le t-1$, the graph $G_d$ is obtained by taking the union of the two graphs $B_{d -t/2 + 1}$ and $H$.
\end{prf}

We are now ready to present a necessary condition for recovery.
In particular, the next theorem
implies that the sufficient condition presented in Theorem~\ref{th: IP} is tight.

\begin{theorem}
\label{th: IP tight}
Let $n$ be a positive integer divisible by eight, and let $\din,\dout$ be nonnegative integers such that $\din \le n/2-1$ and $\dout \le n/2$.
If $\din - \dout \le n/4 - 1$, there is a graph $G$ with $n$ nodes and a planted bisection with parameters $\din,\dout$ for which the minimum bisection problem
does not recover the planted bisection.
\end{theorem}

\begin{prf}
Let $n,\din, \dout$ satisfy the assumptions in the statement.
We explain how to construct a graph $G$ with $n$ nodes, and the planted bisection $V_1, V_2$
with parameters $\din, \dout$.
Furthermore we show that 
the minimum bisection problem over $G$ does not recover the planted bisection.

The notation that we use in this proof is the same used in the
proof of Theorem~\ref{th: IP} and we refer the reader to Fig.~\ref{fig: flip} for an illustration.
In our instances the set $V_1$ is the union of the two disjoint and nonempty sets of nodes $U_1, W_1$.
Symmetrically, the set $V_2$ is the union of the two disjoint sets of nodes $U_2, W_2$.
We will always have $|U_1| = |U_2|$, thus $|W_1| = |W_2|$.
In order to define our instances, we will use the sets of edges $R_1, R_2, R$, $S_1, S_2, S$, and $T_1, T_2, T$, as defined in the proof of Theorem~\ref{th: IP}.


By assumption, we have $\dout \in \{0,\dots,n/2\}$.
We subdivide the proof into two separate cases: $\din \le \dout -1$ and $\din \ge \dout$.
In all cases below, the constructed instance will have $|S| \le |R|$.
Thus, for these instances, a solution of the minimum bisection problem with cost no larger than the planted bisection $V_1, V_2$ is given by the bisection $U_1 \cup W_2, U_2 \cup W_1$.

\emph{Case 1: $\din \le \dout -1$.}
Let $G[V_1]$ and $G[V_2]$ be $\din$-regular graphs.
Since $|V_1|$ and $|V_2|$ are even and $0 \le \din \le n/2 - 1$, these graphs exist due to Lemma~\ref{lem: regular}.
Let $T \cup R$ be the edge set of a bipartite $\dout$-regular graph on nodes $V_1 \cup V_2$ with bipartition $V_1,V_2$.
Note that this graph exists due to Lemma~\ref{lem: regular}, since $0 \le \dout \le n/2$.
Then $G_0$ is $\dout$-regular.
Let $U_1$ contain only one node of $V_1$, and $U_2$ contain only one node of $V_2$.
Hence the sets $W_1, W_2$ contain $n/2-1$ nodes each.
It is then simple to verify that $|S| = 2\din$ and $|R| \ge 2 (\dout -1)$.
Thus we have
$|S| = 2 \din \le 2 (\dout -1) \le |R|$.

\emph{Case 2: $\din \ge \dout$.}
In the instances that we construct in this case we have that all sets $U_1, W_1, U_2, W_2$ have cardinality $n/4$. 
We distinguish between two subcases.
In the first we assume $\din \le n/4 - 1$, while in the second subcase we have $\din \ge n/4$.

\emph{Case 2a: $\din \ge \dout$ and $\din \le n/4 - 1$.}
Let $G[U_1], G[W_1], G[U_2], G[W_2]$ be $\din$-regular graphs.
Since $|U_1| = |W_1| = |U_2| = |W_2|$ is even and $0 \le \din \le n/4 - 1$, these graphs exist due to Lemma~\ref{lem: regular}.
Let $S_1, S_2$ be empty.
Then $G[V_1]$ and $G[V_2]$ are $\din$-regular.

Let $T_1$ be the edge set of a bipartite $\dout$-regular graph on nodes $U_1 \cup U_2$ with bipartition $U_1, U_2$, which exists due to Lemma~\ref{lem: regular}, since $0 \le \dout \le \din \le n/4 - 1$.
Symmetrically, let $T_2$ be the edge set of a bipartite $\dout$-regular graph on nodes $W_1 \cup W_2$ with bipartition $W_1, W_2$.
Furthermore, let $R_1, R_2$ be empty.
Then $G_0$ is $\dout$-regular.

In the instance constructed we have $|S| = |R| = 0$.

\emph{Case 2b: $\din \ge \dout$ and $\din \ge n/4$.}
Let $G[U_1], G[W_1], G[U_2], G[W_2]$ be complete graphs, which are $(n/4-1)$-regular.
Let $S_1$ be the edge set of a bipartite $(\din-n/4+1)$-regular graph on nodes $U_1 \cup W_1$ with bipartition $U_1, W_1$.
This bipartite graph exists due to Lemma~\ref{lem: regular} since
\begin{align*}
1 = \cancel{n/4} - \cancel{n/4} + 1 \le \din - n/4 + 1 \le n/2 - 1 - n/4 + 1 = n/4,
\end{align*}
where the second inequality holds by assumption of the theorem.
Symmetrically, let $S_2$ be the edge set of a bipartite $(\din-n/4+1)$-regular graph on nodes $U_2 \cup W_2$ with bipartition $U_2, W_2$.
Then $G[V_1]$ and $G[V_2]$ are $\din$-regular.

Employing Lemma~\ref{lem: regular} similarly to the previous paragraph, we let $R_1$ be the edge set of a bipartite $(\din-n/4+1)$-regular graph on nodes $U_1 \cup W_2$ with bipartition $U_1, W_2$, and let $R_2$ be the edge set of a bipartite $(\din-n/4+1)$-regular graph on nodes $U_2 \cup W_1$ with bipartition $U_2, W_1$.
Furthermore, let $T_1$ be the edge set of a bipartite $(\dout - \din+n/4-1)$-regular graph on nodes $U_1 \cup U_2$ with bipartition $U_1,U_2$.
This bipartite graph exists due to Lemma~\ref{lem: regular} since
\begin{align*}
0 \le \dout - \din+n/4-1 \le \cancel{\din} - \cancel{\din} + n/4 - 1 = n/4 - 1,
\end{align*}
where the first inequality holds by assumption of the theorem.
Symmetrically, let $T_2$ be the edge set of a bipartite $(\dout - \din+n/4-1)$-regular graph on nodes $W_1 \cup W_2$ with bipartition $W_1, W_2$.
Then $G_0$ is $\dout$-regular.

Note that in the instance defined we have $|S_1| = |S_2| = |R_1| = |R_2|$, thus $|S| = |R|$.
\end{prf}

\sctn{Linear programming relaxation}
\label{sec: lp}

In order to present the LP relaxation for the minimum bisection problem, we first give an equivalent formulation of this problem.
Given a graph $G=(V,E)$ with $n$ nodes, the \emph{cut vector} corresponding to a partition of $V$ into two sets is defined as the vector $x \in \{0,1\}^{\binom{n}{2}}$ with $x_{ij} = 0$ for any $i<j \in V$ with nodes $i,j$ in the same set of the partition, and with $x_{ij} = 1$ for any $i<j \in V$ with nodes $i,j$ in different sets of the partition.
Throughout the paper, we use the notations $x_{ij}$ and $x_{ji}$ interchangeably; that is, we set $x_{ji} := x_{ij}$ for all $i < j \in V$.
We define $a_{ij}=1$ for every $(i,j) \in E$ and $a_{ij}=0$ for every $i,j \in V$ with $(i,j) \notin E$.
Then the minimum bisection problem can be written as the problem of minimizing the linear function $\sum_{(i,j) \in E} x_{ij} = \sum_{1 \le i< j \le n} a_{ij} x_{ij}$ over all cut vectors subject to $\sum_{1\leq i < j\leq n} {x_{ij}} = n^2/4$, where the equality constraint ensures that the two sets in the partition have the same cardinality.

The most well-known LP relaxation of the minimum bisection problem is the so-called~\emph{metric relaxation}, given by
\begin{align}
\label{primal}
\tag{LP-P}
\min_x  \quad & \sum_{1\leq i < j \leq n} {a_{ij} x_{ij}} \nonumber \\
\sub \quad & x_{ij} \leq x_{ik} + x_{jk} && \forall \; {\rm distinct} \; i,j, k \in [n], \; i < j \label{tri1}\\
& x_{ij} + x_{ik} + x_{jk} \leq 2 && \forall \; 1\leq i < j < k \leq n \label{tri2}\\
& \sum_{1 \leq i < j \leq n} {x_{ij}} = n^2/4\label{eqsize}.
\end{align}
Problem~\eqref{primal}
is obtained by first convexifying the feasible region of the minimum bisection problem by requiring $x$ to be in the \emph{cut polytope}, \ie the convex hull of all cut vectors, while satisfying the condition that the two sets in the partition have equal size, enforced by~\eqref{eqsize}.
Subsequently, the cut polytope is outer-approximated by the metric polytope defined by triangle inequalities~\eqref{tri1} and~\eqref{tri2}.
From a complexity viewpoint, Problem~\eqref{primal} can be solved in polynomial time.
From a practical perspective, it is well-understood that by employing some cutting-plane method together with dual Simplex algorithm, Problem~\eqref{primal} can be solved efficiently.

\sbsctn{A sufficient condition for recovery}
\label{sec:LPsuff}

Our goal in this section is to obtain sufficient conditions under which the LP relaxation recovers the planted bisection.
To this end, first, we
obtain conditions under which the cut vector corresponding to the planted bisection is an optimal solution of~\eqref{primal}.
Subsequently, we address the question of uniqueness.
We start by constructing the dual of~\eqref{primal};
we associate variables $\lambda_{ijk}$
with inequalities~\eqref{tri1},
variables $\mu_{ijk}$
with inequalities~\eqref{tri2}, and a variable $\omega$ with the equality constraint~\eqref{eqsize}.
It then follows that the dual of~\eqref{primal} is given by:
\begin{align}
\label{dual}
\tag{LP-D}
\max_{\lambda, \mu, \omega}  \quad & -2\sum_{1\leq i < j < k \leq n} {\mu_{ijk}} - \frac{n^2}{4} \omega\nonumber\\
\sub \quad & a_{ij} + \sum_{k\in [n]\setminus \{i,j\}} (\lambda_{ijk}-\lambda_{ikj}-\lambda_{jki} + \mu_{ijk}) + \omega = 0 && 1 \leq i < j \leq n \label{d1}\\
& \lambda_{ijk} = \lambda_{jik} \geq 0, \quad  \mu_{ijk} = \mu_{ikj} = \mu_{kij} \geq 0 && \forall i \neq j \neq k \in [n], \; i < j   \nonumber
\end{align}
Let $\bar x$ be feasible to the primal \eqref{primal} and $(\bar\lambda, \bar\mu, \bar\omega)$ be feasible to the dual \eqref{dual}.
Then $\bar x$ and $(\bar\lambda, \bar\mu, \bar\omega)$ are optimal,
if they satisfy complementary slackness.
As before let $G$ be a graph with planted bisection $V_1, V_2$. Let $\bar x$ be the cut vector corresponding to the planted bisection.
Without loss of generality, we assume $V_1 = \{1,\dots,\frac  n2\}$, $V_2 = \{\frac  n2 + 1,\dots,n\}$.
Then we have $\bar x_{ij} =0$ for all $i,j \in V_1$ and all $i,j \in V_2$ with $i<j$, and $\bar x_{ij} = 1$ for all $i \in V_1$ and $j \in V_2$.
For notational simplicity, in the remainder of this paper we let $G_1: = G[V_1]$ and $G_2 := G[V_2]$.

\sbsbsctn{Regular graphs}
\label{regg}

In the following, we consider the setting where $G_1$ and $G_2$ are both $\din$-regular graphs for some
$\din \in \{1, \ldots, \frac{n}{2}-1\}$ and that $G_0$ is $\dout$-regular for some $\dout \in \{0,\ldots, \frac{n}{2}\}$,
where we assume $n \geq 4$.
This restrictive assumption significantly simplifies the optimality conditions and enables us to obtain the dual certificate in
closed-form. In the next section, we relax this regularity assumption.

\begin{proposition}
\label{optimal}
Let $G_1$ and $G_2$ be $\din$-regular for some
$\din \in \{1, \ldots, \frac{n}{2}-1\}$ and let $G_0$ be $\dout$-regular for some $\dout \in \{0,\ldots, \frac{n}{2}\}$.
Then the cut vector $\bar x$ corresponding to the planted bisection is an optimal solution of~\eqref{primal} if
\begin{equation}\label{condition}
\din - \dout \geq \frac{n}{4}-1.
\end{equation}
\end{proposition}

\begin{prf}
We prove the statement by constructing a dual feasible point $(\bar\lambda, \bar\mu, \bar\omega)$ that together with $\bar x$ satisfies complementary slackness.
By complementary slackness, at any such point we have $\bar\lambda_{ijk} = 0$ for all $i,j \in V_1$, $k \in V_2$ and for all
$i,j \in V_2$, $k \in V_1$ and $\bar\mu_{ijk} = 0$ for all $i, j, k \in V_1$ and all $i, j, k \in V_2$.
Additionally, we make the following simplifying assumptions:
\begin{itemize}
    \item $\bar\lambda_{ijk} = 0$ for all $i, j, k \in V_1$ and for all $i, j, k \in V_2$ ,
    \item $\bar\lambda_{ikj} = \bar\lambda_{jki} = 0$ for all $i,j \in V_1$ such that $(i,j) \notin E_1$ and for all $k \in V_2$, and $\bar\lambda_{kij} = \bar\lambda_{kji} = 0$ for all $i,j \in V_2$ such that $(i,j) \notin E_2$ and for all $k \in V_1$,
    \item $\bar\mu_{ijk} = 0$ for all $(i,j) \in E_1$ and $k \in V_2$,
and $\bar\mu_{ijk} = 0$ for all $(i,j) \in E_2$ and $k \in V_1$.
\end{itemize}
%
%
It then follows that the linear system~\eqref{d1} simplifies to following set of equalities:
\begin{itemize}

\item [(I)] For any $i< j \in V_1$ and $(i,j) \notin E_1$:
\begin{equation*}
 \sum_{k \in V_2}{\bar\mu_{ijk}}  + \bar\omega = 0.
\end{equation*}

\item [(II)] For any $i<j \in V_2$ and $(i,j) \notin E_2$:
\begin{equation*}
\sum_{k \in V_1}{\bar\mu_{ijk}}  + \bar\omega = 0.
\end{equation*}

\item [(III)] For any $i \in V_1$ and $j \in V_2$ such that $(i,j) \notin E_0$:
\begin{equation*}
  \sum_{\substack{k \in V_1: \\ (i,k) \in E_1}}{(\bar\lambda_{ijk}-\bar\lambda_{kji})}
+ \sum_{\substack{k \in V_2: \\ (j,k) \in E_2}}{(\bar\lambda_{ijk}-\bar\lambda_{ikj})}
+ \sum_{\substack{k \in V_1: \\ (i,k) \notin E_1}} {\bar\mu_{ijk}}
+ \sum_{\substack{k \in V_2: \\ (j,k) \notin E_2}}{\bar\mu_{ijk}}
+  \bar\omega = 0.
\end{equation*}

\item [(IV)] For any $(i,j) \in E_0$:
\begin{equation*}
1 + \sum_{\substack{k \in V_1 : \\ (i,k) \in E_1}}{(\bar\lambda_{ijk}-\bar\lambda_{kji})}
+ \sum_{\substack{k \in V_2 : \\ (j,k) \in E_2}}{(\bar\lambda_{ijk}-\bar\lambda_{ikj})}
+ \sum_{\substack{k \in V_1 : \\ (i,k) \notin E_1}} {\bar\mu_{ijk}}
+ \sum_{\substack{k \in V_2 : \\ (j,k) \notin E_2}}{\bar\mu_{ijk}}
+  \bar\omega = 0.
\end{equation*}

\item [(V)] For any $(i,j) \in E_1$:
\begin{equation*}
1 - \sum_{k \in V_2}{(\bar\lambda_{i k j} + \bar\lambda_{j k i})}  + \bar\omega = 0.
\end{equation*}

\item [(VI)] For any $(i,j) \in E_2$:
\begin{equation*}
1 - \sum_{k \in V_1}{(\bar\lambda_{k i j} + \bar\lambda_{k j i})}  + \bar\omega = 0.
\end{equation*}

\end{itemize}

First, to satisfy condition~(I) (resp. condition~(II)),
for each $i<j \in V_1$ such that $(i,j) \notin E_1$ and $k \in V_2$ (resp. $i<j \in V_2$ such that $(i,j) \notin E_2$ and $k \in V_1$), let
\begin{equation}\label{wtf10}
\bar\mu_{ijk} = -\big(\frac{a_{ik} + a_{jk}}{2}\big) \frac{\bar\omega}{\dout}.
\end{equation}
Substituting~\eqref{wtf10} in condition~(I) yields
$$
-\frac{\bar\omega}{2\dout} \sum_{k \in V_2}(a_{ik} + a_{jk}) + \bar \omega = -\frac{\bar\omega}{2\dout} 2 \dout + \bar \omega = 0,
$$
where the first equality follows from $\dout$-regularity of $G_0$.
It then follows that to satisfy $\bar\mu_{ijk} \geq 0$, we must have
\begin{equation}\label{cond1}
\bar\omega \leq 0.
\end{equation}
We will return to this condition once we determine $\bar \omega$.
Next, to satisfy conditions~(III), for each $(i, j) \in E_1$ and $k \in V_2$, let
\begin{equation}\label{wtf11}
\bar\lambda_{ikj}-\bar\lambda_{jki} = \big(\frac{a_{ik}-a_{jk}}{2}\big) \frac{\bar\omega}{\dout},
\end{equation}
and for each $(i, j) \in E_2$ and $k \in V_1$, let
\begin{equation}\label{wtf11p}
\bar\lambda_{kij}-\bar\lambda_{kji} = \big(\frac{a_{ik}-a_{jk}}{2}\big) \frac{\bar\omega}{\dout}.
\end{equation}
%
Substituting~\eqref{wtf10}, ~\eqref{wtf11}, and~\eqref{wtf11p} in condition~(III),  for each $i \in V_1$ and $j \in V_2$ such that $(i,j) \notin E_0$, we obtain
\begin{align*}
-\frac{\bar \omega}{2\dout}\sum_{k \in V_1}{a_{jk}}
-\frac{\bar \omega}{2\dout} \sum_{k \in V_2}{a_{ik}}
+  \bar\omega = -\frac{\bar\omega}{2} -\frac{\bar\omega}{2} + \bar\omega = 0,
\end{align*}
where the first equality follows from $\dout$-regularity of $G_0$.

Next we choose $\bar \omega$ so that condition~(IV) is satisfied.
Substituting~\eqref{wtf10},\eqref{wtf11}, and ~\eqref{wtf11p} in condition~(IV), for each $(i,j) \in E_0$ we obtain
\begin{align*}
& 1+\frac{\bar \omega}{2 \dout} \Big(\sum_{\substack{k \in V_1,\\ (i,k) \in E_1}}{(1-a_{jk})}+\sum_{\substack{k \in V_2, \\ (j,k) \in E_2}}{(1-a_{ik})}
- \sum_{\substack{k \in V_1, \\ (i,k) \notin E_1}}{(1+a_{jk})}- \sum_{\substack{k \in V_2, \\ (j,k) \notin E_2}}{(1+a_{ik})}\Big)+\bar\omega \\
=&  1+\frac{\bar \omega}{2 \dout} \Big(\sum_{\substack{k \in V_1,\\ (i,k) \in E_1}}{1}-\sum_{\substack{k \in V_1,\\ (i,k) \notin E_1}}{1} - \sum_{\substack{k \in V_1 \setminus \{i\}}}{a_{jk}}+\sum_{\substack{k \in V_2, \\ (j,k) \in E_2}}{1}-\sum_{\substack{k \in V_2, \\ (j,k) \notin E_2}}{1}
- \sum_{\substack{k \in V_2 \setminus \{j\}}}{a_{ik}}\Big)+\bar\omega \\
=&  1+\frac{\bar \omega}{2 \dout} \Big(\din-(\frac{n}{2}-1-\din)-\dout+1+\din-(\frac{n}{2}-1-\din)-\dout+1\Big)+\bar\omega\\
=&1 + \frac{\bar\omega}{\dout}(2 \din + 2-\dout-\frac{n}{2}) +\bar\omega =0,
\end{align*}
where the second equality follows from $\din$-regularity of $G_1, G_2$
and $\dout$-regularity of $G_0$. Hence, to satisfy condition~(IV) we must have:
\begin{equation}\label{omegaeq}
\bar\omega = \frac{2\dout}{n-4\din-4}.
\end{equation}
It then follows that to satisfy condition~\eqref{cond1}
we must have $\din \geq \frac{n}{4} -1$, which is clearly implied by~\eqref{condition}.

To complete the proof, by symmetry, it suffices to find nonnegative $\bar\lambda$ satisfying condition~(V) together with
relation~\eqref{wtf11}. For each $(i,j) \in E_1$ and for each $k \in V_2$, letting
\begin{align}\label{add1}
&\bar \lambda_{ikj} = \max \Big\{0, \; \big(\frac{a_{ik}-a_{jk}}{2}\big)\frac{\bar \omega}{\dout}\Big\} + \gamma_{ij}, \quad \bar \lambda_{jki} = \max \Big\{0, \; \big(\frac{a_{jk}-a_{ik}}{2}\big)\frac{\bar \omega}{\dout}\Big\} + \gamma_{ij},
\end{align}
for some $\gamma_{ij} \geq 0$,
it follows that condition~(V) can be satisfied if
$$
-\frac{\bar\omega}{2\dout}\sum_{k\in V_2}{|a_{ik}-a_{jk}|} \leq 1+ \bar \omega.
$$
By $\dout$-regularity of $G_0$, we have $\sum_{k\in V_2}{|a_{ik}-a_{jk}|} \leq 2 \dout$ for all $i,j \in V_1$; hence the above inequality holds if
\begin{equation} \label{last}
1 + 2 \bar\omega \geq 0.
\end{equation}
By~\eqref{omegaeq}, inequality~\eqref{last} is equivalent to condition~\eqref{condition} and this completes the proof.
\end{prf}

We now provide a sufficient condition under which the cut vector corresponding to the planted bisection is the~\emph{unique} optimal solution of~\eqref{primal}.
In~\cite{Man79}, Mangasarian studies necessary and sufficient conditions for a solution of an LP to be unique.
Essentially, he shows that an LP solution is unique if and only if it remains a
solution to each LP obtained by an arbitrary but sufficiently
small perturbation of its objective function. Subsequently, he gives a number of equivalent characterizations, one of which, stated
below, turned out to be useful in our context.

\begin{theorem} [Theorem~2, part~(iv) in~\cite{Man79}]\label{tool}
Consider the linear program:
\begin{eqnarray*}
& \min_x \quad  & p^T x \\
& \sub & Ax= b, \quad Cx \geq d,
\end{eqnarray*}
where $p$, $b$ and $d$ are vectors in $\R^n$, $\R^m$ and $\R^k$ respectively, and $A$ and $C$ are $m \times n$ and $k \times n$ matrices respectively. The dual of this linear program is given by:
\begin{eqnarray*}
& \max_{u,v} \quad & b^T u + d^T v \\
& \sub & A^T u + C^T v= p, \quad v \geq 0.
\end{eqnarray*}
Let $\bar x$ be an optimal solution of the primal and let $(\bar u, \bar v)$ be an optimal solution of the dual.
Let $C_i$ denote the $i$th row of $C$. Define $K= \{ i: C_i \bar x= d_i, \bar v_i >0 \}$ and $L= \{ i: C_i \bar x= d_i, \bar v_i = 0 \}$.
Let $C_K$ and $C_L$, be matrices whose rows are $C_i$, $i \in K$ and $i \in L$,
respectively. Then $\bar x$ is a unique optimal solution if and only if there exist no $x$ different from the zero vector
satisfying
\begin{equation*}
A x = 0, \quad C_K x = 0, \quad C_L x \geq 0.
\end{equation*}
\end{theorem}

Using the above characterization, we now present a sufficient condition for uniqueness of the optimal solution of~\eqref{primal}.

\begin{proposition}
\label{uniqueness}
Let $G_1$ and $G_2$ be $\din$-regular and let
$G_0$ be $\dout$-regular. Then the LP relaxation recovers the planted bisection, if
\begin{equation}\label{cond3}
\din - \dout \geq \frac{n}{4}.
\end{equation}
\end{proposition}

\begin{prf}
We start by characterizing the index sets $K$ and $L$ defined in the statement of Theorem~\ref{tool} for Problem~\eqref{primal}.
By~\eqref{wtf10}, $\bar \mu_{ijk} > 0$ for the following set of triplets $(i,j,k)$:
\begin{enumerate}
\item all $i , j \in V_1$ such that $(i,j) \notin E_1$ and for all $k \in V_2$ such that $(i,k) \in E_0$ or $(j,k) \in E_0$
\item all $i , j \in V_2$ such that $(i,j) \notin E_2$ and for all $k \in V_1$ such that $(i,k) \in E_0$ or  $(j,k) \in E_0$,
\end{enumerate}
and $\bar \mu_{ijk} = 0$, otherwise. Moreover, by condition~\eqref{cond3}, we have $\gamma_{ij} > 0$ which in turn implies $\bar \lambda_{ikj}, \bar \lambda_{jki} > 0$ for all $(i,j) \in E_1$, $k \in V_2$ and by symmetry $\bar \lambda_{k i j}, \bar\lambda_{k j i} > 0$ for all $(i,j) \in E_2$, $k \in V_1$ and $\bar \lambda$s equal zero, otherwise.

Hence, by Theorem~\ref{tool} and the proof of Proposition~\ref{optimal} it suffices to show that
there exists no $x$ different from the zero vector satisfying the following
\begin{itemize}
\item [(i)] $\sum_{1 \leq i < j \leq n} {x_{ij}} = 0$.

\item [(ii)] For each $(i, j) \in E_1$ and each $k \in V_2$ (symmetrically for each $(i,j) \in E_2$ and each $k \in V_1$):
             $x_{ik} = x_{ij} + x_{jk}$, $x_{jk} = x_{ij} + x_{ik}$ and $x_{ij} + x_{jk} + x_{ik} \leq 0$.

\item [(iii)] For each $(i, j) \notin E_1$ and each $k \in V_2$ such that $(i,k) \in E_0$ or $(j,k) \in E_0$
              (symmetrically, for each $(i, j) \notin E_2$ and each $k \in V_1$ such that $(i,k) \in E_0$ or $(j,k) \in E_0$):
              $x_{ij} + x_{jk} + x_{ik} = 0$, $x_{ik} \leq x_{ij} + x_{jk}$, $x_{jk} \leq x_{ij} + x_{ik}$.

\end{itemize}

To obtain a contradiction, assume that there exists a nonzero $x$ satisfying conditions~(i)-(iii).
From condition~(ii) it follows that
\begin{equation}\label{u1}
x_{ij} = 0, \quad \forall (i,j) \in E_1 \cup E_2,
\end{equation}
and for each $(i,j) \in E_1$  (resp. $(i,j) \in E_2$) we have $x_{ik} = x_{jk} \leq 0$
for all $k \in V_2$ (resp. $k \in V_1)$.    Recall that a Hamiltonian cycle is a cycle that visits each node exactly once and a graph is called Hamiltonian if it has a Hamiltonian cycle.
Dirac (1952) proved that a simple graph with $m$ nodes ($m \geq 3$) is Hamiltonian if every node has degree $\frac{m}{2}$ or greater.
By assumption~\eqref{cond3}, we have $\din \geq \frac{n}{4}$. Since $G_1$ and $G_2$ are $\din$-regular graphs with $\frac{n}{2}$ nodes, by Dirac's result, they are both Hamiltonian. Now consider a Hamiltonian cycle in $G_1$ (resp. $G_2$)
denoted by $v_1 v_2 \ldots v_{n/2} v_1$. It then follows that for each $k \in V_2$  we have
$x_{v_i k} = x_{v_j k}$ for all $i,j \in \{1,\ldots, n/2\}$. Symmetrically, for each $k \in V_1$,
$x_{u_i k} = x_{u_j k}$ for all $i,j \in \{1,\ldots, n/2\}$, where $u_1 u_2 \ldots u_{n/2} u_1$
denotes a Hamiltonian cycle in $G_2$.
Consequently, we have
\begin{equation}\label{u2}
x_{ij} = \alpha \leq 0, \quad \forall i \in V_1, \; j \in V_2.
\end{equation}

From condition~(iii), for each
$(i,j) \notin E_1$ and each $k \in V_2$ such that $(i,k) \in E_0$ or $(j,k) \in E_0$, we have
$x_{ij} = -x_{ik} - x_{jk} = -2 \alpha$, where the second equality follows from~\eqref{u2}.
By $\dout$-regularity of $G_0$ and also by symmetry, we conclude that
\begin{equation}\label{u3}
x_{ij} = -2 \alpha, \quad \forall (i,j): i < j \in V_1, \; (i,j) \notin E_1, \; {\rm or} \; i < j \in V_2, \; (i,j) \notin E_2
\end{equation}
Finally substituting~\eqref{u1},~\eqref{u2}, and~\eqref{u3} into condition~(i), we obtain
\begin{align*}
\Big(\frac{n}{2}\Big)^2 \alpha - \Big(\frac{n}{2}\Big)\Big(\frac{n}{2}-\din-1\Big) (2 \alpha) = \alpha n \Big(\din + 1 -\frac{n}{4}\Big) = 0,
\end{align*}
which by assumption~\eqref{cond3} holds only if $\alpha = 0$. However this implies that $x$ has to be the zero vector, which contradicts our assumption. Hence, under assumption~\eqref{cond3}, the LP recovers the planted bisection.
\end{prf}

The results of Proposition~\ref{uniqueness} and Theorem~\ref{th: IP tight} indicate that the worst-case recovery condition for the LP relaxation is~\emph{tight}, provided that the input graph is characterized in terms of $n,\din,\dout$, $G_1$ and $G_2$ are $\din$-regular, and $G_0$ is $\dout$-regular.
Indeed, this is a surprising result as the triangle inequalities constitute a very small fraction of the facet-defining inequalities
for the cut polytope.

In~\cite{BoeDenStr20} the authors present a sufficient condition
under which the minimum ratio-cut problem recovers a planted partition, in terms of the spectrum of the adjacency matrix of $G$.
It can be shown that a similar recovery condition can also be obtained for the minimum bisection problem. Such a condition yields stronger recovery
guarantees than that of Theorem~\ref{th: IP} for sparse random graphs.
At the time of this writing, we are not able to construct a dual certificate using the spectrum of the adjacency matrix, and this remains an interesting subject for future research.

\sbsbsctn{General graphs}

We now relax the regularity assumptions of the previous section.
To this end, we make use of the next two lemmata. In the following we characterize a graph $G$ by the three
subgraphs $(G_0, G_1, G_2)$, as defined before.

\begin{lemma}\label{l1}
Suppose that the LP relaxation recovers the planted bisection for $(G_0, G_1, G_2)$. Then it also recovers
the planted bisection for $(\tilde G_0, G_1, G_2)$ where $E(\tilde G_0) \subset E(G_0)$.
\end{lemma}

\begin{prf}
Let $|E(G_0) \setminus E(\tilde G_0)| = m$ for some $m \geq 1$ and denote by $\bar x$ the cut vector corresponding to
the planted bisection of $(G_0, G_1, G_2)$.
Since the LP relaxation recovers the planted bisection for $(G_0, G_1, G_2)$, the optimal value of this LP is equal to $|E(G_0)|$.
Moreover, the objective value of the LP relaxation for $(\tilde G_0, G_1, G_2)$ at $\bar x$
is equal to $|E(G_0)| -m$. If $\bar x$ is not uniquely optimal for the latter LP, there exists a different solution $\tilde x$ that gives an objective value of $|E(G_0)| -m - \delta$ for some $\delta \geq 0$. Now let us compute the objective value of the LP relaxation
for $(G_0, G_1, G_2)$ at $\tilde x$; we get
$\sum_{e \in E(G_1)} {\tilde x_e}+ \sum_{e \in E(G_2)} {\tilde x_e} + \sum_{e \in E(\tilde G_0)} {\tilde x_e} +
 \sum_{e \in E(G_0) \setminus E(\tilde G_0)} {\tilde x_e} =|E(G_0)| -m - \delta + \sum_{e \in E(G_0) \setminus E(\tilde G_0)} {\tilde x_e}$. Since
 the planted bisection is recovered for $(G_0, G_1, G_2)$ we must have $|E(G_0)| -m - \delta + \sum_{e \in E(G_0) \setminus E(\tilde G_0)} {\tilde x_e} > |E(G_0)|$; \ie $\sum_{e \in E(G_0) \setminus E(\tilde G_0)} {\tilde x_e} > m + \delta$. However this is not possible since we have $x_{ij} \leq 1$ for all $i,j$.
\end{prf}

\begin{lemma}\label{l2}
Suppose that the LP relaxation recovers the planted bisection for $(G_0, G_1, G_2)$. Then it also recovers
the planted bisection for $(G_0, \tilde G_1, \tilde G_2)$ where $E(\tilde G_1) \supseteq E(G_1)$ and $E(\tilde G_2) \supseteq E(G_2)$.
\end{lemma}

\begin{prf}
Denote by $\bar x$ the cut vector corresponding to the planted bisection of $(G_0, G_1, G_2)$.
Since the LP relaxation recovers the planted bisection for $(G_0, G_1, G_2)$,
the optimal value of the LP is equal to $|E(G_0)|$.
It then follows that $\bar x$ gives an objective value of $|E(G_0)|$ for the LP corresponding to $(G_0, \tilde G_1, \tilde G_2)$ as well. If
$\bar x$ is not uniquely optimal for the latter LP, it means there exists a different solution $\tilde x$ that gives an objective value of $|E(G_0)| - \delta$ for some $\delta \geq 0$. Now compute the objective value of the LP for
$(G_0, G_1, G_2)$ at $\tilde x$; we obtain
$\sum_{e \in E(G_1)} {\tilde x_e}+ \sum_{e \in E(G_2)} {\tilde x_e} + \sum_{e \in E(G_0)} {\tilde x_e} =
\sum_{e \in E(\tilde G_1)} {\tilde x_e}- \sum_{e \in E(\tilde G_1) \setminus E(G_1)} {\tilde x_e}+ \sum_{e \in E(\tilde G_2)} {\tilde x_e}
- \sum_{e \in E(\tilde G_2) \setminus E(G_2)} {\tilde x_e}
+ \sum_{e \in E(G_0)} {\tilde x_e} =
|E(G_0)| - \delta - \sum_{e \in E(\tilde G_1) \setminus E(G_1)} {\tilde x_e} -\sum_{e \in E(\tilde G_2) \setminus E(G_2)} {\tilde x_e}$. Since the planted bisection is recovered for $(G_0, G_1, G_2)$ we must have $|E(G_0)| - \delta - \sum_{e \in E(\tilde G_1) \setminus E(G_1)} {\tilde x_e} -\sum_{e \in E(\tilde G_2) \setminus E(G_2)} {\tilde x_e} > |E(G_0)|$; \ie
$\sum_{e \in E(\tilde G_1) \setminus E(G_1)} {\tilde x_e} +\sum_{e \in E(\tilde G_2) \setminus E(G_2)} {\tilde x_e} < - \delta$. However this is not possible since we have $x_{ij} \geq 0$ for all $i,j$.
\end{prf}

Utilizing Lemma~\ref{l1}, Lemma~\ref{l2} and Proposition~\ref{uniqueness}, we now present a sufficient condition
for recovery of the planted bisection in general graphs. As in Section~\ref{sec: bisection}, we denote by
$\din$ the minimum node degree of $G_1 \cup G_2$ and we denote by $\dout$ the maximum node degree of $G_0$.

\begin{theorem}
\label{deter}
Assume that $G_1$ and $G_2$ contain $\din$-regular subgraphs on the same node set
and assume that $G_0$ is a subgraph of a $\dout$-regular bipartite graph with the same bipartition.
Then the LP relaxation recovers the planted bisection, if
$\din - \dout \geq \frac{n}{4}$.
\end{theorem}

We should remark that in general, the assumptions of Theorem~\ref{deter} are restrictive;
consider a graph $G_1$ with the minimum node degree $\din$. Then it is simple to construct instances
which do not contain a $\din$-regular subgraph on the same node set. Similarly, in general, a bipartite graph
$G_0$ with maximum node degree $\dout$ is not a subgraph of a $\dout$-regular bipartite graph with the same bipartition.
However, as we show in Section~\ref{sec: SBM}, for certain random graphs, these assumptions hold with high probability.

\sbsctn{A necessary condition for recovery}
\label{sec:LPnec}

Consider Problem~\eqref{primal}; in this section, we present a necessary condition for the recovery of the planted bisection using this LP relaxation.
To this end, we construct a feasible point of the LP whose corresponding objective value, under certain conditions,
is strictly smaller than that of the planted
bisection. This in turn implies that the LP relaxation does not recover the planted bisection.

Before proceeding further,
we introduce some notation
that we will use to present our result.
Consider a graph $G = (V, E)$, and let $\rho_G: V^2 \to \mathbb{Z}_{\geq 0} \cup \{+\infty\}$ denote the pairwise distances in $G$; \ie $\rho_G(i, j)$ is the length of the shortest path from $i$ to $j$ in $G$, with $\rho_G(i, j) = +\infty$ if $i$ and $j$ belong to different connected components of $G$.
Define the \emph{diameter} and \emph{average distance} in $G$ as follows:
\begin{align*}
  \rho_{\max}(G) &\colonequals \max_{i, j \in V} \rho_G(i, j), \\
  \rho_{\tavg}(G) &\colonequals \avg_{i, j \sim \mathrm{Unif}(V)} [\rho_G(i, j)] = \frac{2}{n^2}\sum_{1 \leq i < j \leq n} \rho_G(i, j).
\end{align*}

Let us also give some intuition behind our argument, which will also give some explanation for the relevance of these quantities.
Problem~\eqref{primal} optimizes, as its name suggests, over a choice of a \emph{metric} on the vertex set of $G$, and seeks to make the distances under this metric small between pairs of points that are adjacent in $G$.
To find a non-integral feasible point with a large value, it is therefore natural to start with the metric $\rho_G$ itself, and try to renormalize it to build a feasible $x$ for Problem~\eqref{primal}.
This is precisely what we will do in the proof of the next result; $\rho_{\tavg}(G)$ will be involved when we normalize $x$ to have a fixed value of $\sum_{i < j} x_{ij}$, and $\rho_{\max}(G)$ when we adjust our construction to satisfy the other ``triangle inequalities'' included in the problem.

\begin{theorem}
    \label{lp-nonrecovery-cond}
    Let $G = (V, E)$ be a connected graph on $n \geq 5$ nodes.
    Define
    \begin{equation}\label{cG}
    c(G) \colonequals \max\left\{0, \frac{3\rho_{\max}(G) - 4\rho_{\tavg}(G)}{1 - \frac{4}{n}}\right\}.
    \end{equation}
    If
    \begin{equation}
      \label{eq:lp-nonrecovery-cond}
      \frac{1 + c(G)}{2\rho_{\tavg}(G) + 2c(G)(1 - \frac{1}{n})} |E| < |E_0|,
  \end{equation}
  then LP relaxation~\eqref{primal} does not recover the planted bisection.
\end{theorem}

\begin{prf}
Define the point $\tilde{x} \in \mathbb{R}^{\binom{n}{2}}$ as
$$\tilde{x}_{ij} = \frac{\rho_G(i, j) + c(G)}{2\rho_{\tavg}(G) + 2c(G)(1 - \frac{1}{n})} = \frac{\rho_G(i, j) + c(G)}{\frac{4}{n^2}\sum_{1 \leq i < j \leq n}(\rho_G(i, j) + c(G))}.$$
We first show that  $\tilde{x}$ is feasible for~\eqref{primal}.
By construction $\tilde{x}$ satisfies the equality constraint~\eqref{eqsize}.
Also, since $\rho_G$ is a metric and $c \geq 0$, $\tilde{x}$ satisfies inequalities~\eqref{tri1}.
    Thus it suffices to show that $\tilde{x}$ satisfies inequalities~\eqref{tri2}.
    To do this, we bound:
    \begin{equation}\label{feas}
    \tilde{x}_{ij} + \tilde{x}_{ik} + \tilde{x}_{jk} \leq \frac{3}{2} \cdot \frac{\rho_{\max}(G) + c(G)}{\rho_{\tavg}(G)+ c(G)(1 - \frac{1}{n})}. \end{equation}
    Two cases arise:
    \begin{itemize}
    \item [(i)] $c(G) = 0$:
    in this case, we have $\rho_{\max}(G) \leq \frac{4}{3}\rho_{\tavg}(G)$ and hence we can further bound the right-hand side of inequality~\eqref{feas} as:
    \[ \tilde{x}_{ij} + \tilde{x}_{ik} + \tilde{x}_{jk} \leq \frac{3}{2} \cdot \frac{\rho_{\max}(G)}{\rho_{\tavg}(G)} \leq \frac{3}{2} \cdot \frac{4}{3} = 2. \]
    \item [(ii)] $c(G) > 0$:
    in this case, we have $\rho_{\max}(G) > \frac{4}{3} \rho_{\tavg}(G)$ and $c(G) = (3\rho_{\max}(G) - 4\rho_{\tavg}(G)) / (1 - \frac{4}{n})$;
    hence we can further bound the right-hand side of inequality~\eqref{feas} as:
    \begin{align*}
      \tilde{x}_{ij} + \tilde{x}_{ik} + \tilde{x}_{jk}
      &\leq \frac{3}{2} \cdot \frac{\rho_{\max}(G)(1 - \frac{4}{n}) + 3\rho_{\max}(G) - 4\rho_{\tavg}(G)}{\rho_{\tavg}(G)(1 - \frac{4}{n}) + (3\rho_{\max}(G) - 4\rho_{\tavg}(G))(1 - \frac{1}{n})} \\
      &= \frac{3}{2} \cdot \frac{4\rho_{\max}(G)(1 - \frac{1}{n}) - 4\rho_{\tavg}(G)}{3\rho_{\max}(G)(1 - \frac{1}{n})-3\rho_{\tavg}(G)} \\
      &= 2.
    \end{align*}
    \end{itemize}
    Thus in either case $\tilde{x}$ satisfies inequalities~\eqref{tri2}, implying $\tilde x$ is feasible.
    It then follows that the objective value of~\eqref{primal} at $\tilde x$, given by $\frac{1 + c(G)}{\rho_{\tavg}(G) + c(G)(1 - \frac{1}{n})} \cdot \frac{|E|}{2}$,
     provides an upper bound on its optimal value.
    Moreover, the objective value of~\eqref{primal} corresponding to the planted bisection is equal to $|E_0|$.
    Hence, if condition~\eqref{eq:lp-nonrecovery-cond} holds,
    the LP relaxation does not recover the planted bisection.
\end{prf}

The no-recovery condition~\eqref{eq:lp-nonrecovery-cond} depends on $c(G)$ which in turn depends on the relative values of
$\rho_{\tavg}(G)$ and $\rho_{\max}(G)$.
Clearly, for any graph $G$ we have $\rho_{\tavg}(G) \leq \rho_{\max}(G)$.
As we detail in Section~\ref{sec: SBM}, for random graphs, very often these two quantities are quite close together, whereby $|E_0|$ must be quite small compared to $|E|$ in order for recovery to be possible.
In short, this result lets us infer that, since typical distances between nodes in random graphs are comparable to the maximum such distance, the LP seldom succeeds in recovering the planted partition.

Before proceeding to those arguments, to illustrate the basic idea we explain how this analysis looks for a specific type of deterministic graph.
Let $G = (V, E)$ be a graph such that $\rho_G(i, j) \leq 2$ for every $i, j \in V$, $\rho_G(i, j) = 2$ for some $i, j \in V$.
That is, $G$ is not the complete graph, but any pair of non-adjacent nodes in $G$ have a common neighbor.
In this case, we may explicitly compute $\rho_{\max}(G) = 2$ and
\[ \rho_{\tavg}(G) = \frac{2}{n^2}\left(|E| + 2\left(\binom{n}{2} - |E|\right)\right) = 2\left(1 - \frac{1}{n}\right) - \frac{2}{n^2}|E|. \]
Thus $3\rho_{\max}(G) - 4\rho_{\tavg}(G) = -2 + \frac{8}{n} + \frac{8}{n^2}|E|$.
If moreover $n$ is sufficiently large and 
$|E| \leq \frac{n^2}{4} - n$,
we have $c(G) = 0$.
Thus the non-recovery condition \eqref{eq:lp-nonrecovery-cond} reduces to $|E_0| / |E| > \frac{1}{2\rho_{\tavg}(G)}$, where $\rho_{\tavg}(G) \in (1, 2)$.
Therefore, in this type of graph our analysis implies that if a sufficiently large constant fraction of edges cross the planted bisection, then the LP cannot recover the planted bisection.

\sctn{The stochastic block model}
\label{sec: SBM}

The \emph{stochastic block model (SBM)} is a probability distribution over instances of the planted bisection problem, which in recent years has been studied extensively as a benchmark for community detection (see~\cite{Abbe18} for an extensive survey).
The SBM can be seen as an extension of the \emph{Erd\H{o}s-R\'{e}nyi (ER)} model for random graphs.
In the ER model, edges are placed between each pair of nodes independently with probability $p$.
In the SBM, there is an additional community structure built in to the random graph.
In this paper, we focus on the simplest case, the \emph{symmetric SBM on two communities}.
In this model, for an even number $n$ of nodes, we first choose a bisection $V_1, V_2$ of the set of nodes uniformly at random.
We then draw edges between pairs of nodes $i, j \in V_1$ or $i, j \in V_2$ with probability $p$, and draw edges between pairs of nodes $i \in V_1$ and $j \in V_2$ with probability $q$.
As in the ER model, all edges are chosen independently.
In addition, we always assume $p>q$, so that the average connectivity within each individual community is stronger than that between different communities.


We denote by $\G_{n, p}$ the ER model on $n$ nodes with edge probability $p$, and by $\G_{n, p, q}$ the SBM on $n$ nodes  with edge probabilities $p$ and $q$.
We write $G \sim \G_{n, p}$ or $G \sim \G_{n, p, q}$ for $G$ drawn from either distribution.
Generally, we consider probability parameters depending on $n$, $p = p(n)$ and $q = q(n)$, but for the sake of brevity we often omit this dependence. When $p$ and $q$ are in fact constants independent of $n$, we refer to the graph $G \sim \G_{n,p,q}$ as \emph{very dense}.
When instead $p$ and $q$ scale as
\begin{equation}\label{DenseG}
p = \alpha n^{-\omega}, \qquad
q = \beta n^{-\omega},
\end{equation}
where $\alpha$, $\beta$ and $\omega$ are parameters independent of $n$ such that $\alpha > \beta > 0$ and $0 <\omega < 1$, we refer to the graph $G \sim \G_{n,p,q}$ as \emph{dense} and finally when $p$ and $q$ scale as
\begin{equation}\label{SparseG}
p = \alpha \frac{\log n}{n}, \qquad
q = \beta \frac{\log n}{n},
\end{equation}
where $\alpha$ and $\beta$ are parameters independent of $n$ such that $\alpha > \beta > 0$, we refer to the graph $G \sim \G_{n,p,q}$ as having \emph{logarithmic degree}.

Let us briefly review some of the existing literature on recovery in the SBM, to establish a baseline for our results.
In~\cite{AbbBanHal16}, the authors show that the minimum bisection problem recovers the planted bisection with high probability for
very dense and dense graphs as long as $p > q$. For graphs of logarithmic degree, they show that the minimum bisection problem recovers the planted bisection with high probability if and only if $\sqrt{\alpha} - \sqrt{\beta} > \sqrt{2}$.
This gives an ``information-theoretic limit'' on recovery, in the sense that if $\sqrt{\alpha} - \sqrt{\beta} < \sqrt{2}$, then no algorithm (efficient or not) can recover the planted bisection with high probability.
Moreover, the authors consider an SDP relaxation of the minimum bisection problem, and conjecture that in the logarithmic regime the SDP relaxation recovers the planted bisection
whenever the minimum bisection problem does.
The authors of \cite{AbbBanHal16} provided a proof of a partial result in this direction, and the subsequent works \cite{HajWuXu16,Ban18} gave complete proofs.
In Section~\ref{sdpTim}
we show that the deterministic recovery condition in \cite{Ban18} implies that the SDP relaxation also recovers the planted bisection for very dense and dense graphs up to the information theoretic threshold.
The final picture emerging from this analysis is that the SDP relaxation is ``as good as'' the information-theoretically optimal integer program.

In this section, we investigate the theoretical performance of the LP relaxation for recovery in the SBM.
We show that in any of the three regimes described above, the LP relaxation does~\emph{not} recover the planted bisection with high probability
in the full regions of parameters where it is information-theoretically possible to do so.
More specifically, we first obtain a sufficient condition under which the LP relaxation
recovers the planted bisection for very dense graphs.
Subsequently, we obtain non-recovery conditions for very dense and dense graphs.
Finally we show that for graphs of logarithmic degree, the LP relaxation fails to recover the planted bisection with high probability.
In short, the LP relaxation is not ``as good as'' the SDP relaxation
in recovering the communities in the SBM.

\sbsctn{Recovery in SBM}
\label{sec:RecoverSBM}
In this section, we consider very dense random graphs 
and obtain a sufficient condition under which the LP relaxation recovers the planted bisection with high probability.
To this end, we first show that in the very dense regime, the regularity assumptions of Theorem~\ref{deter} are essentially not restrictive.
Subsequently, we employ Theorem~\ref{deter} to obtain a recovery guarantee for the very dense SBM.
The following propositions provide sufficient conditions under which bipartite ER graphs and general ER graphs
are contained in or contain regular graphs of comparable average degrees.
In the remainder of the paper, for a sequence of events $A_n$ depending on $n$, we say $A_n$ occurs \emph{with high probability} if $\prob[A_n] \to 1$ as $n \to \infty$.

\begin{proposition}\label{prop cool}
There exists a constant $C_{\reg, 2} > 0$ such that, for any $q \in (0, 1)$, with high probability $G \sim \G_{n, 0, q}$ is a subgraph of a $d_{\reg}$-regular bipartite graph
with the same bipartition, where $d_{\reg} = \frac{qn}{2} + C_{\reg, 2}\sqrt{\frac n2 \log \frac n2}$.
\end{proposition}

\begin{proposition}\label{prop cooler}
There exists a constant $C_{\reg, 1} > 0$ such that, for any $p \in (0, 1)$, with high probability $G \sim \G_{n, p}$ contains as a subgraph a $d_{\reg}$-regular graph on the same node set, where $d_{\reg} = pn - C_{\reg, 1}\sqrt{n\log n}$.
\end{proposition}

We note that $\frac{qn}{2}$ and $pn$ are the expected degrees of any given vertex in $G \sim \mathcal{G}_{n, 0, q}$ and $G \sim \mathcal{G}_{n, p}$, respectively.
The proofs of Propositions~\ref{prop cool} and \ref{prop cooler} are given in
Sections~\ref{sec: regularbi} and \ref{sec: regularg},
respectively, and consist of performing a probabilistic analysis of conditions for the existence of ``graph factors'' as developed in classical results of graph theory literature \cite{tutte52}.
As such, these results are of independent interest to both optimization and applied probability communities.

\smallskip

By combining Theorem~\ref{deter} and Propositions~\ref{prop cooler} and~\ref{prop cool},  we obtain
a sufficient condition for recovery in the very dense regime.

\begin{theorem}\label{SBMdense}
    Let $0 < q < p < 1$ with $p - q > \frac{1}{2}$.
    Then, the LP relaxation recovers the planted bisection in $G \sim \G_{n, p, q}$ with high probability.
\end{theorem}

\begin{prf}
As before, denote by $G_1$ and $G_2$ the subgraphs of $G$ induced by $V_1$ and $V_2$, respectively,
and denote by $G_0$ the bipartite subgraph of $G$ of all edges between $V_1$ and $V_2$.
We start by regularizing the bipartite graph $G_0$.
By Proposition~\ref{prop cool}, with high probability, it is possible to add edges to $G_0$ to make it $\dout$-regular, for $\dout = \frac{q n}{2} + C_{\reg, 2} \sqrt{\frac{n}{2} \log (\frac{n}{2})}$ and $C_{\reg, 2}$ a universal constant.
Similarly, by Proposition~\ref{prop cooler}, with high probability, it is possible to remove edges from $G_1$ and $G_2$ to make them $\din$-regular, for $\din = \frac{pn}{2} - C_{\reg, 1} \sqrt{\frac{n}{2} \log (\frac{n}{2})}$ and $C_{\reg, 1}$ another universal constant.
From Theorem~\ref{deter} it follows that the LP relaxation recovers the planted bisection with high probability, if
$\din - \dout \geq \frac{n}{4}$. Substituting for the values of $\din$ and $\dout$ specified above gives that, with high probability, it suffices to have
\begin{align*}
\frac{pn}{2} -\frac{q n}{2} - (C_{\reg, 1} + C_{\reg, 2}) \sqrt{\frac{n}{2} \log \left(\frac{n}{2}\right)}\geq \frac{n}{4}.
\end{align*}
Dividing both side of the above inequality by $n$ and letting $n \rightarrow \infty$ then gives the result, since $p - q > \frac{1}{2}$ by assumption.
\end{prf}

\begin{figure}
    \centering
    \includegraphics[scale=0.35, trim=20mm 0mm 140mm 0mm, clip]{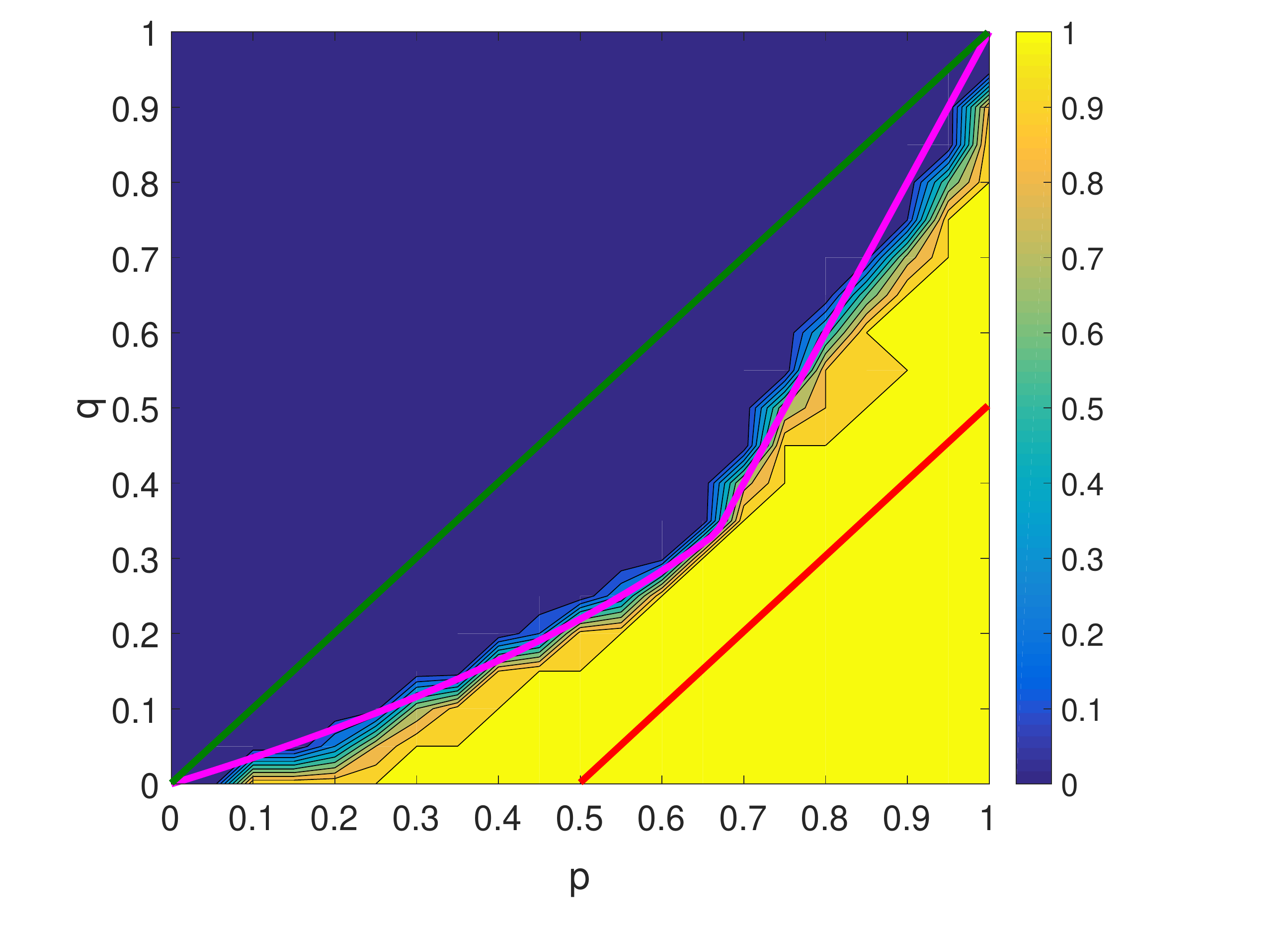}
 \caption{The empirical probability of success of the
LP relaxation in recovering the planted bisection in the very dense regime.
For each fixed pair $(p,q)$, we fix $n = 100$ and the number of trials to be $20$. Then, for each fixed $(p,q)$, we count the number of times the LP relaxation recovers the planted bisection.
Dividing by the number of trials, we obtain the empirical probability
of success. In red, we plot the threshold for recovery of the LP as given by Theorem~\ref{SBMdense}.
In magenta, we plot the threshold for failure of the LP relaxation as given by Part~1 of Theorem~\ref{noRecoveryDense}.
In green, we plot the information-theoretic limit for recovery in the very dense regime.
}
\label{figure2}
\end{figure}

To better understand the tightness of the sufficient condition of Theorem~\ref{SBMdense}, we conduct a simulation, depicted in Figure~\ref{figure2}. As can be seen in this figure, as an artifact of our proof technique, the condition $p - q > \frac{1}{2}$
is somewhat suboptimal.
While a better dual certificate construction scheme may improve this condition, as we can see from this figure, such an improvement is not going to be significant.

\sbsctn{Non-recovery in SBM}
\label{sec:NonRecoverSBM}

We now consider random graphs 
in various 
regimes
and provide conditions under which the recovery in these regimes
is not possible. Our basic technique is to evaluate both sides of the condition~\eqref{eq:lp-nonrecovery-cond} in Theorem~\ref{lp-nonrecovery-cond} asymptotically for the SBM, and to verify that in suitable regimes, the condition holds with high probability.
To this end, for a graph $G \sim \mathcal{G}_{n, p, q}$, we need to bound $c(G)$ defined by~\eqref{cG}.
This in turn amounts to bounding the quantities $\rho_{\max}(G)$ and $\rho_{\tavg}(G)$
in different regimes.
Luckily the distance distributions of ER graphs have been studied extensively in the literature.
Utilizing the existing work on ER graphs, we obtain similar results for the SBM.
We refer the reader to Section~\ref{distances}
for statements and proofs.

Our first result gives a collection of non-recovery conditions for very dense and dense graphs.

\begin{theorem}\label{noRecoveryDense}
    Let $p = \alpha n^{-\omega}$ and $q = \beta n^{-\omega}$ for $\alpha, \beta > 0$ and $\omega \in [0, 1)$, with $\alpha, \beta \in (0, 1)$ if $\omega = 0$.
    Suppose that one of the following conditions holds:
    \begin{enumerate}
    \item $\omega = 0$ and
        \begin{equation}\label{noR1}
            \beta > \max\left\{\frac{3 - \alpha - \sqrt{(3 - \alpha)^2 - 4\alpha}}{2}, \; 2\alpha - 1\right\}.
        \end{equation}
    \item $\omega \in (0, 1)$ with $\frac{1}{1 - \omega} \notin \mathbb{N}$ and
        \begin{equation}\label{noR2}
            \beta > \frac{1}{2\lceil \frac{1}{1 - \omega} \rceil - 1} \cdot \alpha.
        \end{equation}
        \item $\omega = \frac{1}{2}$ and
        \begin{equation}\label{noR3}
            \beta > \max\left\{ \frac{1 - 2e^{-\alpha^2}}{2 - e^{-\alpha^2}},\; \frac{1}{3 + 2e^{-\alpha^2}}\right\} \cdot \alpha.
        \end{equation}
    \item $\omega \in (0, 1)$ with $\frac{1}{1 - \omega} \in \mathbb{N} \setminus \{2\}$ and
        \begin{equation}\label{noR4}
            \beta > \frac{1}{2 \cdot \frac{1}{1 - \omega} - 1 + 2\exp(-\alpha^{\frac{1}{1 - \omega}})} \cdot \alpha.
        \end{equation}
    \end{enumerate}
    Then, with high probability, the LP fails to recover the planted bisection.
\end{theorem}

\begin{prf}
First, we bound the fraction of edges across the bisection in the SBM.
Let $p, q = \Omega(\log n / n)$ and having $\lim_{n \to \infty}p / q$ existing and taking a value in $(0, 1]$.
Then, for any $\epsilon > 0$, an application of Hoeffding's inequality to both $|E|$ and $|E_0|$ shows that
\begin{equation} \label{edge-counts-sbm}
    \lim_{n \to \infty}\prob\left[(1 - \epsilon) \frac{q}{p + q} \leq \frac{|E_0|}{|E|} \leq (1 + \epsilon)\frac{q}{p + q}\right] = 1,
\end{equation}
where we note that $\frac{\avg |E_0|}{\avg |E|} = \frac{q}{p + q} + o(1)$, explaining the appearance of this quantity.

Next, define
\begin{equation}
    \label{eq:def-bG}
    b(G) \colonequals \frac{1 + c(G)}{2\rho_{\tavg}(G) + 2c(G)(1 - \frac{1}{n})}.
\end{equation}
By Theorem~\ref{lp-nonrecovery-cond}, if $|E_0| / |E| > b(G)$, the LP does not recover the planted bisection.
Thus our task consists in obtaining an upper bound on $b(G)$. First, notice that $b(G)$ is increasing in $c(G)$ since
$\rho_{\tavg}(G) \leq \frac{2}{n^2} \cdot \binom{n}{2} = 1 - \frac{1}{n}$. Moreover, $b(G)$ is decreasing in $\rho_{\tavg}(G)$. Hence, to upper bound $b(G)$ it suffices to upper bound $\rho_{\max}(G)$ and to lower bound $\rho_{\tavg}(G)$.
We will treat the four cases individually.

\smallskip

    \emph{Condition 1: $\omega = 0$.}
    By Proposition~\ref{prop:dist-dn-sbm-very-dense}, for any $\epsilon > 0$, the distance distribution in $G$ satisfies, with high probability,
    \begin{equation} \label{eq:dist-dn-sbm-very-dense-pf}
        (1 - \epsilon)\left(2 - \frac{p + q}{2}\right) \leq \rho_{\tavg}(G) \leq \rho_{\max}(G) = 2.
    \end{equation}
    On the event that this holds, we have
    \[ 3\rho_{\max}(G) - 4\rho_{\tavg}(G) \leq 6 - 4(1 - \epsilon)\left(2 - \frac{p + q}{2}\right) \leq 2(p + q - 1) + 8\epsilon. \]
    Therefore, we have, for sufficiently small $\epsilon$, with high probability,
    \[ c(G) \leq \frac{1}{1 - \frac{4}{n}}\left(2\max\{0, p + q - 1\} + 8\epsilon\right) \leq 2\max\{0, p + q - 1\} + 9\epsilon. \]
    We now consider two cases depending on the output of the maximum in the first term.

    \emph{Case 1.1: $p + q \leq 1$.} In this case, with high probability $c(G) \leq 9\epsilon$.
    By \eqref{eq:dist-dn-sbm-very-dense-pf}, we also have, with high probability, $\rho_{\tavg}(G) \geq (1 - \epsilon) \frac{4 - p - q}{2}$.
    Thus from the definition of $b(G)$, for $\epsilon$ sufficiently small, with high probability,
    \[ b(G) \leq (1 + O(\epsilon)) \frac{1}{4 - p - q}. \]
    Therefore, by~\eqref{eq:lp-nonrecovery-cond} and~\eqref{edge-counts-sbm} and by taking $\epsilon$ sufficiently small,
    it suffices to have
    \[ \frac{1}{4 - p - q} < \frac{q}{p + q}.\]
    The above is equivalent to $q > \frac{1}{2}(3 - p - \sqrt{(3 - p)^2 - 4p})$ upon solving the inequality for $q$, completing the proof.

    \emph{Case 1.2: $p + q > 1$.} In this case, with high probability $c(G) \leq 2(p + q - 1) + 9\epsilon$.
    For $\epsilon$ sufficiently small, with high probability,
    \[ b(G) \leq (1 + O(\epsilon)) \frac{1 + 2(p + q - 1)}{4 - p - q + 4(p + q - 1)} = (1 + O(\epsilon)) \frac{2p + 2q - 1}{3p + 3q}. \]
    As in Case 1.1, it then suffices to have
    \[ \frac{2p + 2q - 1}{3p + 3q} < \frac{q}{p + q}. \]
    The above is equivalent to $q > 2p - 1$ upon solving the inequality for $q$, completing the proof.

\smallskip

\emph{Condition 2: $\omega \in (0, 1)$ with $\frac{1}{1 - \omega} \notin \mathbb{N}$.}
    By Proposition~\ref{prop:dist-dn-sbm-dense}, in this regime for any $\epsilon > 0$, the distance distribution in $G$ with high probability satisfies
    \[ (1 - \epsilon)\left\lceil\frac{1}{1 - \omega}\right\rceil \leq \rho_{\tavg}(G) \leq \rho_{\max}(G) = \left\lceil\frac{1}{1 - \omega}\right\rceil. \]
    On this event, we have $3\rho_{\max}(G) - 4\rho_{\tavg}(G) \leq -\lceil \frac{1}{1 - \omega} \rceil + O(\epsilon)$, so for sufficiently small $\epsilon$ this quantity is negative, whereby $c(G) = 0$.
    In this case, we have
    \[ b(G) = \frac{1}{2\rho_{\tavg}(G)} \leq (1 + O(\epsilon))\frac{1}{2\lceil \frac{1}{1 - \omega} \rceil}, \]
    and thus it suffices to have
    \[ \frac{\beta}{\alpha + \beta} > \frac{1}{2\lceil \frac{1}{1 - \omega} \rceil}. \]
    The above condition is equivalent to condition~\eqref{noR2} after solving the inequality for $\beta$, completing the proof.

\smallskip

    \emph{Condition 3: $\omega = \frac{1}{2}$.}
    In this case, $\frac{1}{1 - \omega} = 2$.
    By Proposition~\ref{prop:dist-dn-sbm-exc-dense}, in this regime for any $\epsilon > 0$, the distance distribution in $G$ with high probability satisfies
    \[ (1 - \epsilon)(2 + \exp(-\alpha^{2})) \leq \rho_{\tavg}(G) \leq \rho_{\max}(G) = 3. \]
    On this event, we have
    $$
      3\rho_{\max}(G) - 4\rho_{\tavg}(G) \leq 9 - 4(1 - \epsilon)(2 + \exp(-\alpha^{2})) \leq 1 - 4\exp(-\alpha^2) + O(\epsilon).
    $$
    We then consider two cases:

    \emph{Case 3.1: $\exp(-\alpha^2) \geq \frac{1}{4}$.}
    In this case, we have $3\rho_{\max}(G) - 4\rho_{\tavg}(G) = O(\epsilon)$, whereby $c(G) = O(\epsilon)$ with high probability as well.
    Thus,
    \[ b(G) \leq (1 + O(\epsilon)) \frac{1}{2\rho_{\tavg}(G)} \leq (1 + O(\epsilon))\frac{1}{4 + 2\exp(-\alpha^2)}. \]
    It then suffices to have
    \[ \frac{\beta}{\alpha + \beta} > \frac{1}{4 + 2\exp(-\alpha^2)}, \]
    which, upon solving the inequality for $\beta$ is equivalent to $\beta > \frac{1}{3 + 2\exp(-\alpha^2)} \cdot \alpha$.

    \emph{Case 3.2: $\exp(-\alpha^2) < \frac{1}{4}$.}
    In this case, we have $c(G) \leq (1 + O(\epsilon))(1 - 4\exp(-\alpha^2))$ with high probability.
    Thus,
    \[ b(G) \leq (1 + O(\epsilon))\frac{1 + (1 - 4\exp(-\alpha^2))}{2\rho_{\tavg}(G) + 2(1 - 4\exp(-\alpha^2))} \leq (1 + O(\epsilon))\frac{1 - 2\exp(-\alpha^2)}{3 - 3\exp(-\alpha^2)}. \]
    It then suffices to have
    \[ \frac{\beta}{\alpha + \beta} > \frac{1 - 2\exp(-\alpha^2)}{3 - 3\exp(-\alpha^2)}, \]
    which, upon solving the inequality for $\beta$ is equivalent to $\beta > \frac{1 - 2\exp(-\alpha^2)}{2 - \exp(-\alpha^2)} \cdot \alpha$.

\smallskip

    \emph{Condition 4: $\omega \in (0, 1)$ with $\frac{1}{1 - \omega} \in \mathbb{N} \setminus \{2\}$.}
    In this case, $\omega = \frac{k - 1}{k}$ for some $k \in \mathbb{N}$ with $k \geq 3$, and $\frac{1}{1 - \omega} = k$.
    By Proposition~\ref{prop:dist-dn-sbm-exc-dense}, in this regime for any $\epsilon > 0$, the distance distribution in $G$ with high probability satisfies
    \[ (1 - \epsilon)(k + \exp(-\alpha^{k})) \leq \rho_{\tavg}(G) \leq \rho_{\max}(G) = k + 1. \]
    On this event, we have
    $$
      3\rho_{\max}(G) - 4\rho_{\tavg}(G)
      \leq 3 - (1 - 4\epsilon)k - 4(1 - \epsilon)\exp(-\alpha^{k})
      \leq 12\epsilon - 4(1 - \epsilon)\exp(-\alpha^{k}),
    $$
    where the second inequality is valid since $k \geq 3$. It then follows that
    for $\epsilon$ sufficiently small, $3\rho_{\max}(G) - 4\rho_{\tavg}(G) < 0$, and hence for such $\epsilon$ with high probability $c(G) = 0$.
    On this event, we have
    \[ b(G) = \frac{1}{2\rho_{\tavg}(G)} \leq (1 + O(\epsilon))\frac{1}{2k + 2\exp(-\alpha^k)}, \]
    and thus it suffices to check that
    \[ \frac{\beta}{\alpha + \beta} > \frac{1}{2k + 2\exp(-\alpha^k)}. \]
    The above inequality is equivalent to inequality~\eqref{noR4} after solving the inequality for $\beta$, completing the proof.
\end{prf}

\begin{figure}
\label{fig:lp-negative-dense}
    \begin{center}
        \includegraphics[scale = 0.6, trim=0 13 0 0]{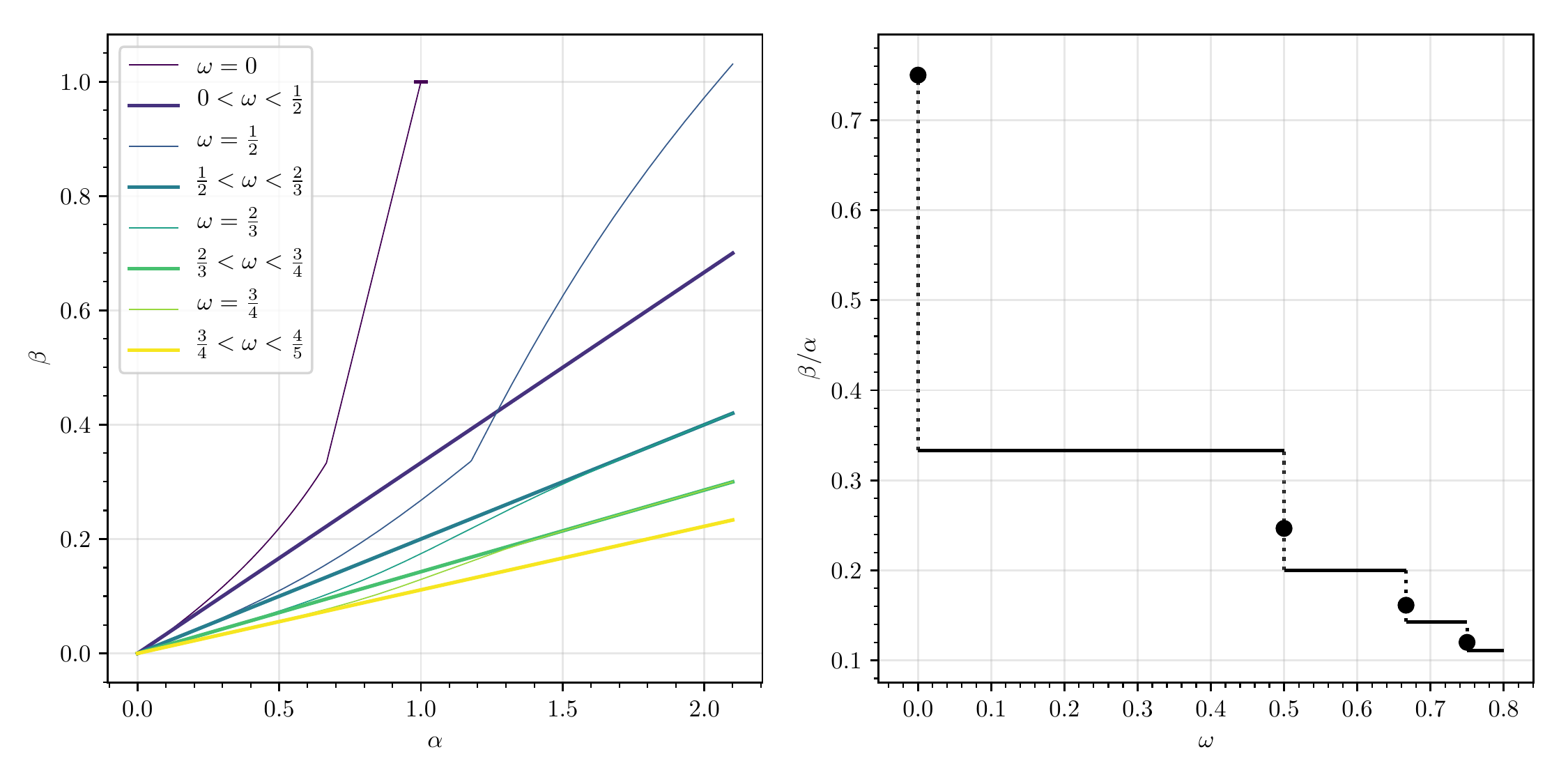} 
    \end{center}
\caption{Non-recovery thresholds as given by Theorem~\ref{noRecoveryDense} for SBM in the very dense and dense regimes. For each choice of $\omega$, our results prove that with high probability the LP does not recover the planted bisection for all choices of $(\alpha, \beta)$ lying above the given curve in the left panel. The curve for $\omega = 0$ is identical to that compared with numerical results in Figure~\ref{figure2}. In the right panel, we show the value of $\beta / \alpha$ at the non-recovery threshold value of $\beta$ for a fixed $\alpha = $ 0.8 as $\omega$ varies. In contrast, the same ratio for the SDP relaxation is always equal to 1.}
\end{figure}

The non-recovery thresholds of Theorem~\ref{noRecoveryDense} are plotted in Figure~\ref{fig:lp-negative-dense} for several values of $\omega$.
Moreover, condition~\eqref{noR1} is plotted in Figure~\ref{figure2} as well. As can be seen from this figure, our non-recovery threshold
in the very dense regime is almost in perfect agreement with our empirical observations.

\smallskip

Finally, we show that the LP relaxation fails to recover the planted bisection in the logarithmic regime.

\begin{theorem}\label{noRecoverySparse}
    Suppose $p = \alpha \frac{\log n}{n}$ and $q = \beta \frac{\log n}{n}$ with $\alpha, \beta > 0$.
    Then, with high probability, the LP relaxation fails to recover the planted bisection in $G \sim \mathcal{G}_{n, p, q}$.
\end{theorem}
\begin{prf}
    If $\frac{\alpha + \beta}{2} \leq 1$, then the statement follows from the information-theoretic bound of \cite{AbbBanHal16}.
    Let us suppose that $\frac{\alpha + \beta}{2} > 1$;
    in this case, since $\frac{q}{p + q} = \frac{\beta}{\alpha + \beta} > 0$ is a constant, it suffices to show that for any $\delta > 0$, with high probability, $b(G) < \delta$.
    By Proposition~\ref{prop:dist-dn-sbm-sparse}, in the logarithmic regime, for any $\epsilon > 0$, the distance distribution in $G$ satisfies, with high probability,
    \begin{equation}
        (1 - \epsilon)\frac{\log n}{\log\log n} \leq \rho_{\tavg}(G) \leq \rho_{\max}(G) \leq (1 + \epsilon)\frac{\log n}{\log\log n}.
    \end{equation}
    It then follows that with high probability $c(G) = 0$ and $b(G) \leq \frac{1}{2(1-\epsilon)} \frac{\log \log n}{\log n}$ and hence the result follows.
\end{prf}

  In~\cite{LeiRao99}, the authors consider the metric relaxation of the sparsest cut problem.
  They show that, for constant degree expander graphs, which includes random regular graphs of constant degree, the LP relaxation has an integrality gap of $\Theta(\log n)$.
  We outline this argument: a $d$-regular expander satisfies, by definition, $|E(U, U^c)| \gtrsim d \min\{|U|, |U^c|\}$ for all $U \subseteq V$, with ``$\gtrsim$'' hiding a constant not depending on $n$.
  Thus the sparsest cut objective function is $\min_U \frac{|E(U, U^c)|}{|U| \cdot |U^c|} \gtrsim \frac{d}{n}$.
  On the other hand, the authors of \cite{LeiRao99} show that the value of the metric relaxation is $\lesssim \frac{\log d}{\log n} \cdot \frac{d}{n}$.
  For $d$ constant, this shows the integrality gap of $\Theta(\log n)$.

  For random regular and ER graphs of average degree $d = d(n) \gtrsim \log n$, the same lower and upper bounds for the sparsest cut and its metric relaxation are valid. Up to an adjustment of constants, we also expect them to hold for SBM graphs, since these are nothing but ``non-uniform ER graphs''.
  For $d \sim \log n$, we have $\frac{\log d}{\log n} \sim \frac{\log \log n}{\log n} = o(1)$, so we expect the metric relaxation to have a diverging integrality gap over all SBM graphs.
  Thus it seems plausible that the proof techniques of~\cite{LeiRao99} could reproduce the non-recovery result of Theorem~\ref{noRecoverySparse}.
  On the other hand, for $d \sim n^{1 - \omega}$ as in our dense regimes, we have $\frac{\log d}{\log n} \sim 1 - \omega$ is a constant.
  Since we expect that working with SBM graphs rather than ER graphs also creates a constant order adjustment in the integrality gap, the argument of~\cite{LeiRao99} does not show non-recovery, and the more refined analysis of Theorem~\ref{noRecoveryDense} is necessary.

Our results do still leave open the question of whether the LP relaxation fails to recover the planted bisection in the entirety of the dense regimes $p, q \sim n^{-\omega}$ with $\omega \in (0, 1)$ of the SBM, or whether they undergo a transition from recovery to non-recovery as in the very dense regime $\omega = 0$.
As we detailed at the end of
Section~\ref{regg}, it may be possible to strengthen our deterministic recovery condition in Proposition~\ref{optimal} by characterizing the input graph in terms of the spectrum of its adjacency matrix rather than the two parameters $\din, \dout$.
This could imply that recovery occurs with high probability in part of the dense regime; a full characterization of the regime over which the LP relaxation undergoes a transition from recovery to non-recovery is a subject for future research. Finally we must point out that
triangle inequalities constitute a very small fraction of facet-defining inequalities for the cut polytope and many more classes of facet-defining inequalities are known for this polytope~\cite{dl:97}. It would be interesting to investigate the recovery properties of stronger LP relaxations for the min-bisection problem.


\sctn{Technical proofs}
\label{sec:proofs}

\sbsctn{Proof of Proposition~\ref{prop cool}}
\label{sec: regularbi}

To prove the statement,
we make use of two results.
The first result provides a necessary and sufficient condition
under which a bipartite graph has a $b$-factor.
Consider a graph $G = (V, E)$; given a vector $b \in \Z^V_{\geq 0}$, a \emph{$b$-factor}
is a subset $F$ of $E$ such that $b_v$ edges in $F$ are incident to $v$
for all $v \in V$.
For notational simplicity, given a vector $b \in \Z^V_{\geq 0}$ and a subset $U \subseteq V$, we define $b(U) := \sum_{v \in U} {b_v}$.
Moreover, for $U, U^{\prime} \subseteq V$, we denote by $N(U,U^{\prime})$ the number of edges in $G$
with one endpoint in $U$ and one endpoint in $U^{\prime}$.

\begin{proposition}[Corollary 21.4a in~\cite{SchBookCO}]\label{prop:bfactor-bipartite}
Let $G = (V,E)$ be a bipartite graph and let $b \in \Z^V_{\geq 0}$.
Then $G$ has a $b$-factor if and only if, for each $U \subseteq V$,
\begin{equation*}
    N(U, U) \geq b(U) - \frac{1}{2} b(V).
\end{equation*}
\end{proposition}

The second result that we need is concerned with an upper bound on the maximum node degree in a bipartite dense random graph.

\begin{proposition}
    \label{prop:max-degree-dense-bipartite}
    Let $G \sim \G_{2m, 0, q}$ for some $q \in (0, 1)$.
    Denote by $d_{\max}$ the maximum node degree of $G$.
    Then, there exists a constant $C > 0$ such that, with high probability, $d_{\max} \leq qm + C \sqrt{m\log m}$.
\end{proposition}

\begin{prf}
We introduce the bipartite adjacency matrix of $G$:
let $X \in \R^{V_1 \times V_2}$, where when $i \in V_1$ and $j \in V_2$, then $X_{ij} = 1$ if $i$ and $j$ are adjacent in $G$ and $X_{ij} = 0$ otherwise.
Then, when $G \sim \G_{2m, 0, q}$, $X_{ij}$ has i.i.d.\ entries equal to $1$ with probability $q$ and 0 otherwise.
In particular, $\avg [X_{ij}] = q$.

The degree of $i \in V_1$ is given by $d(i) = \sum_{j \in V_2} X_{ij}$.
In particular, the $d(i)$ are themselves i.i.d.\ random variables.
Since $X_{ij} \in [0, 1]$, Hoeffding's inequality
applies to each $d(i)$, giving, for $C > 0$ a large constant to be fixed later,
\begin{align*}
 \prob\left[d(i) > qm + C \sqrt{m\log m}\right]
      &= \prob\left[\sum_{j \in V_2} (X_{ij} - \avg [X_{ij}]) > C \sqrt{m\log m}\right] \\
      &\leq \exp\left(\frac{-2C^2m\log m}{m}\right)
      = m^{-2C^2}.
\end{align*}
Since each $d(i)$ is independent for distinct $i \in V_1$, we moreover have, for some fixed $i_1 \in V_1$,
    \begin{align*}
      \prob\Big[ d(i) > qm & + C \sqrt{m\log m} \text{ for some } i \in V_1\Big] \\
      &= 1 - \prob\left[ d(i) \leq qm + C \sqrt{m\log m} \text{ for all } i \in V_1\right] \\
      &= 1 - \left(1 - \prob\left[d(i_1) > qm + C \sqrt{m\log m}\right]\right)^m \\
      &\leq 1 - \left(1 - m^{-2C^2}\right)^m
      = m^{-2C^2}\left(\sum_{k = 0}^{m - 1}\left(1 - m^{-2C^2}\right)^k\right) \\
      &\leq m \cdot m^{-2C^2}
      = m^{1 - 2C^2}.
    \end{align*}
Symmetrically, the same bound applies to $V_1$ replaced with $V_2$, so we find
    \begin{equation*}
        \prob\left[d_{\max} > qm + C \sqrt{m\log m}\right] \leq 2m^{1 - 2C^2},
    \end{equation*}
and setting $C$ large enough gives the result.
\end{prf}

\smallskip

We can now proceed with the proof of Proposition~\ref{prop cool}.
We first set some useful notation.
    Let $V$ be the node set of $G$ and let $V_1, V_2$ be the bipartition of $V$, so that $|V_1| = |V_2| = m$.
    For $U \subseteq V$, let $U_1 = U \cap V_1$ and $U_2 = U \cap V_2$, whereby $U$ is the disjoint union of $U_1$ and $U_2$.
    For $v \in V$, let $d(v)$ be the degree of $v$.
    Let $\bar{G}$ be the bipartite graph complement of $G$, with the same node partition.
    Finally, for $U, U^{\prime} \subseteq V$, let $\bar{N}(U, U^{\prime})$ be the number of edges in $\bar{G}$ with one endpoint in $U$ and one endpoint in $U^{\prime}$.

    The statement of the theorem is equivalent to $\bar{G}$ having a $b$-factor for
    \begin{equation*}
        b(v) = d_{\reg} - d(v).
    \end{equation*}

    Let $d_{\max}$ be the maximum node degree of $G$.
    By Proposition~\ref{prop:max-degree-dense-bipartite}, for sufficiently large $C_{\reg, 2}$, we will have $d_{\reg} \geq d_{\max}$ with high probability, and thus $b(v) \geq 0$ with high probability.
    Let $A_{\deg}$ be the event that this occurs.

    On the event $A_{\deg}$, by Proposition~\ref{prop:bfactor-bipartite}, such a $b$-factor exists if and only if, for each $U \subseteq V$,
    \begin{equation}
        \label{eq:bipartite-b-factor-cond}
        \bar{N}(U, U) \geq b(U) - \frac{1}{2}b(V).
    \end{equation}
    Let $A_{U}$ be the event that this occurs.
    Then, whenever the event $A = A_{\deg} \cap \bigcap_{U \subseteq V} A_U$ occurs, $G$ is a subgraph of a $d_{\reg}$-regular bipartite graph.
    Thus it suffices to show that $\prob[A^c] \to 0$ as $n \to \infty$.
    Taking a union bound,
    \begin{equation}
        \label{eq:bipartite-ub-initial}
        \prob[A^c] \leq \prob[A_{\deg}^c] + \sum_{U \subseteq V} \prob[A_U^c].
    \end{equation}
    We have already observed that for $C_{\reg, 2}$ large enough, the first summand tends to zero, so it suffices to control the remaining sum.

    Let us first rewrite each side of the equation \eqref{eq:bipartite-b-factor-cond}.
    For any $U \subseteq V$, we have
    \[ b(U) = d_{\reg}|U| - \sum_{u \in U} d(u) = d_{\reg}|U| - 2 N(U, U) - N(U, U^c), \]
    and therefore
    \[ b(U) - \frac{1}{2}b(V) = d_{\reg}\left(|U| - m\right) - N(U, U) + N(U^c, U^c). \]
    Also,
    \[ \bar{N}(U, U) = |U_1| \cdot |U_2| - N(U, U). \]
    Thus,
    \begin{equation*}
      A_U = \bigg\{\bar{N}(U, U) \geq b(U) - \frac{1}{2}b(V)\bigg\} = \bigg\{ N(U^c, U^c) \leq |U_1| \cdot |U_2| + d_{\reg}\left(m - |U|\right)\bigg\}.
    \end{equation*}
    Note that
    \[ N(U^c, U^c) \leq |V_1 \setminus U_1| \cdot |V_2 \setminus U_2| = \left(m - |U_1|\right)\left(m - |U_2|\right) = |U_1| \cdot |U_2| - m\left(|U| - m\right). \]
    Since, for sufficiently large $m$, $d_{\reg} \leq m$, if $|U| \geq m$ then $A_U$ always occurs, so in fact for such $m$ we have $\prob[A_U^c] = 0$ whenever $|U| \geq m$.
    Thus,
    \begin{equation}
        \label{eq:bipartite-ub-restricted}
        \lim_{m \to \infty} \sum_{\substack{U \subseteq V}}\prob[A_U^c] = \lim_{m \to \infty}\sum_{\substack{U \subseteq V \\ |U| < m}} \prob[A_U^c].
    \end{equation}

    Recall that we set $X \in \R^{V_1 \times V_2}$ to be the bipartite adjacency matrix of $G$.
    When $G \sim \G_{n, 0, q}$, $X_{ij}$ has i.i.d.\ entries equal to $1$ with probability $q$ and 0 otherwise.
    In particular, $\avg [X_{ij}] = q$.
    Also, we may compute edge counts as $N(U, U) = \sum_{i \in U_1}\sum_{j \in U_2} X_{ij}$.

    Therefore, we may rewrite
    \begin{align*}
      A_U
      &= \left\{\sum_{i \in V_1 \setminus U_1}\sum_{j \in V_2 \setminus U_2} X_{ij} \leq |U_1| \cdot |U_2| + d_{\reg}(m - |U|)\right\} \\
      &= \left\{\sum_{i \in V_1 \setminus U_1}\sum_{j \in V_2 \setminus U_2} (X_{ij} - \avg [X_{ij}]) \leq |U_1| \cdot |U_2| - q(m - |U_1|)(m - |U_2|) + d_{\reg}(m - |U|)\right\} \\
      &= \left\{\sum_{i \in V_1 \setminus U_1}\sum_{j \in V_2 \setminus U_2} (X_{ij} - \avg [X_{ij}]) \leq (1 - q)|U_1| \cdot |U_2| + (d_{\reg} - qm)(m - |U|)\right\} \\
      &= \left\{\sum_{i \in V_1 \setminus U_1}\sum_{j \in V_2 \setminus U_2} (X_{ij} - \avg [X_{ij}]) \leq (1 - q)|U_1| \cdot |U_2| + C_{\reg, 2} \sqrt{m\log m}(m - |U|)\right\}.
    \end{align*}
    In this form, since $X_{ij} \in [0, 1]$, Hoeffding's inequality
    (Theorem 2.2.6 in~\cite{VerBookHDP})
    applies and gives
    \begin{align*}
      \prob[A_U^c]
      &= \prob\left[\sum_{i \in V_1 \setminus U_1}\sum_{j \in V_2 \setminus U_2} (X_{ij} - \avg [X_{ij}]) > (1 - q)|U_1| \cdot |U_2| +
      C_{\reg, 2}\sqrt{m\log m}(m - |U|)\right] \\
      &\leq \exp\left( -\frac{2\big((1 - q)|U_1| \cdot |U_2| + C_{\reg, 2}\sqrt{m\log m}(m - |U|)\big)^2}{(m - |U_1|)(m - |U_2|)}\right)
    \end{align*}

    Note further that this bound, for a given $U$, depends only on $|U_1|$ and $|U_2|$.
    Thus, in the sum appearing in \eqref{eq:bipartite-ub-restricted}, we may group terms according to new scalar variables $a = |U_1|$ and $b = |U_2|$ and bound
    \begin{align}
      \sum_{\substack{U \subseteq V \\ |U| < m}} \prob[A_U^c]
      &= \sum_{\substack{a, b \in [m] \\ a + b < m}}\sum_{\substack{U_1 \in \binom{V_1}{a} \\ U_2 \in \binom{V_2}{b}}} \prob[A_{U_1 \cup U_2}^c] \nonumber \\
      &\leq \sum_{\substack{a, b \in [m] \\ a + b < m}}\sum_{\substack{U_1 \in \binom{V_1}{a} \\ U_2 \in \binom{V_2}{b}}} \exp\left( -\frac{2\big((1 - q)ab + C_{\reg, 2}\sqrt{m\log m}(m - a - b)\big)^2}{(m - a)(m - b)}\right) \nonumber \\
      &= \sum_{\substack{a, b \in [m] \\ a + b < m}}\binom{m}{a} \binom{m}{b}\exp\left( -\frac{2\big((1 - q)ab + C_{\reg, 2}\sqrt{m\log m}(m - a - b)\big)^2}{(m - a)(m - b)}\right).
      \label{eq:bipartite-ub-scalar}
    \end{align}

    Now, note that since $(m - a)(m - b) = ab + m(m - a - b)$, for any choice of $a, b \in [m]$ with $a + b < m$ we have either $ab \geq \frac{1}{2}(m - a)(m - b)$ or $m(m - a - b) \geq \frac{1}{2}(m - a)(m - b)$.
    Note also that the latter is equivalent to $m - a - b \geq \frac{1}{2m}(m - a)(m - b)$.
    Notating this decomposition more formally, define
    \begin{align*}
      \I &= \{(a, b) \in [m]^2: a + b < m\} \\
      \I_1 &= \left\{(a, b) \in [m]^2: a + b < m, ab \geq \frac{1}{2}(m - a)(m - b)\right\} \\
      \I_2 &= \left\{(a, b) \in [m]^2: a + b < m, m - a - b \geq \frac{1}{2m}(m - a)(m - b)\right\}.
    \end{align*}
    Then, the above observation says that $\I = \I_1 \cup \I_2$.

    When $(a, b) \in \I_1$, then we have
    \begin{align}
      \frac{2\big((1 - q)ab + C_{\reg, 2}\sqrt{m\log m}(m - a - b)\big)^2}{(m - a)(m - b)}
      &\geq \frac{2(1 - q)^2(ab)^2}{(m - a)(m - b)} \nonumber \\
      &\geq \frac{1}{2}(1 - q)^2(m - a)(m - b), \label{eq:bipartite-decomp-bd-1}
    \end{align}
    and when $(a, b) \in \I_2$, then we have
    \begin{align}
      \frac{2\big((1 - q)ab + C_{\reg, 2}\sqrt{m\log m}(m - a - b)\big)^2}{(m - a)(m - b)}
      &\geq \frac{2C_{\reg, 2}^2m\log m (m - a - b)^2}{(m - a)(m - b)} \nonumber \\
      &\geq \frac{1}{2}C_{\reg, 2}^2 \frac{\log m}{m}(m - a)(m - b). \label{eq:bipartite-decomp-bd-2}
    \end{align}
    For sufficiently large $m$, the right-hand side of \eqref{eq:bipartite-decomp-bd-2} is always smaller than the right-hand side of \eqref{eq:bipartite-decomp-bd-1} for all $(a, b) \in \I$.
    On the other hand, either \eqref{eq:bipartite-decomp-bd-1} or \eqref{eq:bipartite-decomp-bd-2} holds for any $(a, b) \in \I$, since $\I = \I_1 \cup \I_2$.
    Thus for sufficiently large $m$, for all $(a, b) \in \I$, we have
    \begin{equation*}
        \frac{2\big((1 - q)ab + C_{\reg, 2}\sqrt{m\log m}(m - a - b)\big)^2}{(m - a)(m - b)} \geq \frac{1}{2}C_{\reg, 2}^2 \frac{\log m}{m}(m - a)(m - b).
    \end{equation*}

    Using this in \eqref{eq:bipartite-ub-scalar}, we may bound
    \begin{align}
      \lim_{m \to \infty}\sum_{\substack{U \subseteq V \\ |U| < m}} \prob[A_U^c]
      &\leq \lim_{m \to \infty}\sum_{\substack{a, b \in [m] \\ a + b < m}}\binom{m}{a} \binom{m}{b}\exp\left( -\frac{1}{2}C_{\reg, 2}^2 \frac{\log m}{m}(m - a)(m - b)\right). \nonumber
      \intertext{Whenever $a + b < m$, then either $a < m / 2$ or $b < m / 2$, so symmetrizing the sum we have}
      &\leq \lim_{m \to \infty}2\sum_{\substack{a \in [m / 2] \\ b \in [m] \\ a + b < m}}\binom{m}{a} \binom{m}{b}\exp\left( -\frac{1}{2}
      C_{\reg, 2}^2 \frac{\log m}{m}(m - a)(m - b)\right) \nonumber \\
      &\leq \lim_{m \to \infty}2\sum_{\substack{a \in [m / 2] \\ b \in [m] \\ a + b < m}}\binom{m}{a} \binom{m}{b}\exp\left( -\frac{1}{4}
      C_{\reg, 2}^2 \log m \cdot (m - b)\right). \nonumber
      \intertext{Now, note that $\binom{m}{b} = \binom{m}{m - b} \leq m^{m - b} = \exp(\log m \cdot (m - b))$, and $\binom{m}{a} \leq m^a = \exp(\log m \cdot a)$. Moreover, $a < m - b$ in all terms in the summation, so $\binom{m}{a}\binom{m}{b} \leq \exp(2 \log m \cdot (m - b))$. Therefore, we continue}
      &\leq \lim_{m \to \infty}2\sum_{\substack{a \in [m / 2] \\ b \in [m] \\ a + b < m}}\exp\left( -\left(\frac{1}{4}C_{\reg, 2}^2 - 2\right) \log m \cdot (m - b)\right). \nonumber
      \intertext{Finally, there are at most $m^2$ terms in the summation, and $m - b \geq 1$ in all terms, so we finish by bounding, for
      $C_{\reg, 2}$ sufficiently large,}
      &\leq \lim_{m \to \infty}2m^2\exp\left( -\left(\frac{1}{4}C_{\reg, 2}^2 - 2\right) \log m\right) \nonumber \\
      &\leq \lim_{m \to \infty}2\exp\left( -\left(\frac{1}{4}C_{\reg, 2}^2 - 4\right) \log m\right) \label{eq:bipartite-ub-final-sum-bd}
    \end{align}
    Choosing $C_{\reg, 2}$ again sufficiently large, this limit will equal zero.

    Combining \eqref{eq:bipartite-ub-restricted} and \eqref{eq:bipartite-ub-final-sum-bd}, we have therefore found that, for $C_{\reg, 2}$ sufficiently large,
    \begin{equation*}
        \lim_{m \to \infty}\sum_{U \subseteq V} \prob[A_U^c] = 0.
    \end{equation*}
    Using this in \eqref{eq:bipartite-ub-initial}, we have
    \begin{equation*}
        \lim_{m \to \infty} \prob[A^c] \leq \lim_{m \to \infty} \prob[A_{\deg}^c] + \lim_{m \to \infty}\sum_{U \subseteq V} \prob[A_U^c] = 0,
    \end{equation*}
    giving the result.

\sbsctn{Proof of Proposition~\ref{prop cooler}}

\label{sec: regularg}

In this proof we consider a general ER graph $G \sim \G_{n,p}$ with the minimum node degree denoted by $d_{\min}$
and we show that with high probability $G$ contains as a subgraph, a $d_{\reg}$-regular graph with the same node set,
where $d_{\reg}$ is asymptotically equal to $d_{\min}$.
To this end we make use of the following results concerning the existence of $b$-factors
for general graphs, degree distribution of dense random graphs and node connectivity of dense random graphs.
Let $U\subseteq V$; in the following, we denote by $G-U$ the graph
obtained from $G$ by removing all the nodes in $U$
together with all the edges incident with nodes in $U$.

\begin{proposition}[\cite{tutte52}]\label{prop:bfactor-unipartite}
    Let $G = (V, E)$ be a graph and let $f \in \Z^V_{\geq 0}$.
    For disjoint subsets $S, T \subseteq V$, let $\widetilde{Q}(S, T)$ denote the number of connected components $C$ of $G-(S \cup T)$ such that $f(C) + N(C, T) \equiv 1 \pmod 2$.
    Then, $G$ has an $f$-factor if and only if, for all disjoint subsets $S, T \subseteq V$,
    \begin{equation*}
        f(S) - f(T) + \sum_{v \in T}d(v) - N(S, T) - \widetilde{Q}(S, T) \geq 0.
    \end{equation*}
\end{proposition}

\begin{proposition}[Chapter 3 of \cite{FriKar16}]
    \label{prop:max-degree-dense-unipartite}
    Let $G \sim \G_{n, p}$ for some $p \in (0, 1)$.
    Denote by $d_{\min}$ the minimum node degree of $G$.
    Then, there exists a constant $C > 0$ such that, with high probability, $d_{\min} \geq pn - C \sqrt{n\log n}$.
\end{proposition}

Recall that the node connectivity $\kappa(G)$ of a graph $G$ is the minimum number of nodes whose deletion disconnects it, or the total number of nodes for $G$ a complete graph.
It is well-known that $\kappa(G) \leq d_{\min}$, where as before $d_{\min}$ denotes the minimum node degree in $G$.
Surprisingly, for ER graphs, the two quantities coincide with high probability.

\begin{proposition}[Theorem 1 of \cite{BolTho85}]
    \label{prop:vertex-connectivity-unipartite}
    Let $G \sim \G_{n, p}$ for some $p \in (0, 1)$.
    Then, with high probability, $\kappa(G) = d_{\min}$.
\end{proposition}

\smallskip

We proceed with the proof of Proposition~\ref{prop cooler}.
First, for $A \subseteq V$, let $Q(A)$ denote the number of connected components in $G - A$.
Clearly, $\widetilde{Q}(S, T) \leq Q(S \cup T)$ for any disjoint $S, T \subseteq V$ (where $\widetilde{Q}(S, T)$ is as defined in the statement of Proposition~\ref{prop:bfactor-unipartite}).

We will apply Proposition~\ref{prop:bfactor-unipartite} to the function $f(v) = d_{\reg}$ for all $v \in V$.
Without loss of generality, we may suppose $d_{\reg}$ is even.
In this case, the case $S = T = \emptyset$ is satisfied automatically, so we may suppose $|S| + |T| \geq 1$ in the sequel.
Using the above observation, we find that the desired regular subgraph exists provided the following event occurs:
\begin{align*}
  A &\colonequals \bigcap_{\substack{S, T \subseteq V \\ S \cap T = \emptyset \\ |S| + |T| \geq 1}} A_{S, T}, \qquad \text{where }
  A_{S, T} \colonequals \left\{N(S, T) - \sum_{v \in T} d(v) \leq d_{\reg}(|S| - |T|) - Q(S \cup T)\right\}.
\end{align*}
We have
\begin{equation*}
    \sum_{v \in T} d(v) = 2 N(T, T) + N(S, T) +N(V \setminus (S \cup T), T),
\end{equation*}
whereby we may rewrite
\begin{align*}
  A_{S, T}
  &= \bigg\{-2 N(T, T) - N(V \setminus (S \cup T), T) \leq d_{\reg}(|S| - |T|) - Q(S \cup T)\bigg\} \\
  &= \bigg\{- \sum_{\{i, j\} \in \binom{T}{2}}(2X_{ij} - \avg [2X_{ij}]) - \sum_{\substack{i \in V \setminus (S \cup T) \\ j \in T}}(X_{ij} - \avg [X_{ij}]) \\
  &\hspace{2cm} \leq p|T|(|T| - 1) + p|T|(n - |S| - |T|) + d_{\reg}(|S| - |T|) - Q(S \cup T)\bigg\} \\
  &= \bigg\{- \sum_{\{i, j\} \in \binom{T}{2}}(2X_{ij} - \avg [2X_{ij}]) - \sum_{\substack{i \in V \setminus (S \cup T) \\ j \in T}}(X_{ij} - \avg [X_{ij}]) \\
  &\hspace{2cm} \leq \bigg(p(n - |T|) - C_{\reg, 1}\sqrt{n\log n}\bigg)|S| + (C_{\reg, 1}\sqrt{n\log n} - p)|T| - Q(S \cup T)\bigg\},
\end{align*}
where we have performed the manipulation
\begin{align*}
  p|T|(|T| - 1)& + p|T|(n - |S| - |T|) + d_{\reg}(|S| - |T|) \\
  &= p|T|(n - |S| - 1) + (pn - C_{\reg, 1}\sqrt{n\log n})(|S| - |T|) \\
  &= p|S|n - p|S|\cdot|T| - p|T| + C_{\reg, 1}\sqrt{n\log n}(|T| - |S|)\\
               &= p|S|(n - |T|) + C_{\reg, 1}\sqrt{n\log n}(|T| - |S|) - p|T| \\
  &= \bigg(p(n - |T|) - C_{\reg, 1}\sqrt{n\log n}\bigg)|S| + (C_{\reg, 1}\sqrt{n\log n} - p)|T|.
\end{align*}

Now, we give a lower bound for the right-hand side by considering two cases, depending on the size of $|T|$.
If $|T| \leq n - 2C_{\reg, 1}p^{-1}\sqrt{n\log n}$, then $p(n - |T|) - C_{\reg, 1}\sqrt{n\log n} \geq C_{\reg, 1}\sqrt{n\log n}$, and for sufficiently large $n$ we also have $C_{\reg, 1}\sqrt{n\log n} - p \geq \frac{1}{4}C_{1, \reg}\sqrt{n\log n}$, so in this case
\begin{equation*}
    \bigg(p(n - |T|) - C_{\reg, 1}\sqrt{n\log n}\bigg)|S| + (C_{\reg, 1}\sqrt{n\log n} - p)|T| \geq \frac{1}{4}C_{\reg, 1}\sqrt{n\log n}(|S| + |T|).
\end{equation*}
If, on the other hand, $|T| \geq n - 2C_{\reg, 1}p^{-1}\sqrt{n\log n}$, then, since $S$ and $T$ are disjoint $|S| + |T| \leq n$ whereby $|S| \leq n - |T| \leq 2C_{\reg, 1}p^{-1}\sqrt{n\log n}$.
Therefore, for sufficiently large $n$, we will have
\begin{align*}
  \bigg(p(n - |T|)& - C_{\reg, 1}\sqrt{n\log n}\bigg)|S| + (C_{\reg, 1}\sqrt{n\log n} - p)|T| \\
  &\geq \frac{1}{2}C_{\reg, 1}n\sqrt{n\log n} - C_{\reg, 1}\sqrt{n\log n}|S|
  \geq \frac{1}{2}C_{\reg, 1}n\sqrt{n\log n} - 2C_{\reg, 1}^2p^{-1}n\log n \\
  &\geq \frac{1}{4}C_{\reg, 1}n\sqrt{n\log n}
  \geq \frac{1}{4}C_{\reg, 1}\sqrt{n\log n}(|S| + |T|).
\end{align*}

Thus the same bound is obtained in either case, and we find that, for sufficiently large $n$, for all $S, T \subseteq V$ disjoint, we have
\begin{align*}
  A_{S, T} &\supseteq \bigg\{- 2\sum_{\{i, j\} \in \binom{T}{2}}(X_{ij} - \avg [X_{ij}]) - \sum_{\substack{i \in V \setminus (S \cup T) \\ j \in T}}(X_{ij} - \avg [X_{ij}]) \nonumber \\
  &\hspace{2cm}\leq \frac{1}{4}C_{\reg, 1}\sqrt{n\log n}(|S| + |T|) - Q(S \cup T)\bigg\}. 
\end{align*}

Now, note that we always have $Q(S \cup T) \leq n$.
Moreover, if $\kappa(G)$ is the node connectivity of $G$, and $|S| + |T| < \kappa(G)$, then $Q(S \cup T) = 1$.
And, by Propositions~\ref{prop:max-degree-dense-unipartite} and \ref{prop:vertex-connectivity-unipartite}, with high probability $\kappa = d_{\min} \geq pn - C \sqrt{n\log n}$.
Let $A_{\con}$ denote the event that both the equality $\kappa = d_{\min}$ and the subsequent inequality hold in $G$.
Then, for sufficiently large $n$, on the event $A_{\con}$,
\begin{equation*}
        Q(S \cup T) \leq \left\{\begin{array}{ll} 1 & \text{if } |S| + |T| \leq pn / 2 \\ n & \text{otherwise}\end{array}\right\} \leq \frac{2}{p}(|S| + |T|),
\end{equation*}
where in the last step we use that $S$ and $T$ are not both empty.
Thus, again for sufficiently large $n$,
\begin{align*}
  &A_{S, T} \supseteq A_{\con} \cap \bigg\{- 2\sum_{\{i, j\} \in \binom{T}{2}}(X_{ij} - \avg [X_{ij}]) - \sum_{\substack{i \in V \setminus (S \cup T) \\ j \in T}}(X_{ij} - \avg [X_{ij}]) \nonumber \\
  &\hspace{4cm} \leq \frac{1}{8}C_{\reg, 1}\sqrt{n\log n}(|S| + |T|)\bigg\}.
\end{align*}

Now, returning to the main event $A$, we note two special cases that make the left-hand side of the definition of the event inside the intersection above equal to zero: (1)~if $T = \emptyset$, and (2)~if $|T| = 1$ and $S \cup T = V$.
Thus we may neglect both of these cases, and write, for sufficiently large $n$,
\begin{align*}
  &A \supseteq A_{\con} \cap \bigcap_{\substack{S, T \subseteq V \\ S \cap T = \emptyset \\ |T| \geq 2 \text{ or } |S| + |T| < n}}\bigg\{- 2\sum_{\{i, j\} \in \binom{T}{2}}(X_{ij} - \avg [X_{ij}]) - \sum_{\substack{i \in V \setminus (S \cup T) \\ j \in T}}(X_{ij} - \avg [X_{ij}]) \nonumber \\
  &\hspace{6cm} \leq \frac{1}{8}C_{\reg, 1}\sqrt{n\log n}(|S| + |T|)\bigg\}.
\end{align*}
Applying a union bound and Hoeffding's inequality,
we may then bound
\begin{align*}
  \lim_{n \to \infty}\prob[A^c]
  &\leq \lim_{n \to \infty}\prob[A_{\con}^c] + \lim_{n \to \infty}\sum_{\substack{S, T \subseteq V \\ S \cap T = \emptyset \\ |T| \geq 2 \text{ or } |S| + |T| < n}}\exp\left(-2\frac{\frac{1}{64}C_{\reg, 1}^2n\log n(|S| + |T|)^2}{2|T|(|T| - 1) + |T|(n - |S| - |T|)}\right). \nonumber
  \intertext{
  The first limit is zero by our previous remark.
  In analyzing the remaining exponential terms, we first bound the denominator as $2|T|(|T| - 1) + |T|(n - |S| - |T|) \leq 3n|T|$, continuing}
  &\leq \lim_{n \to \infty} \sum_{\substack{S, T \subseteq V \\ S \cap T = \emptyset \\ |T| \geq 2 \text{ or } |S| + |T| < n}}\exp\left(-\frac{1}{96}C_{\reg, 1}^2\frac{\log n \cdot (|S| + |T|)^2}{|T|}\right) \nonumber
  \intertext{and then observe that $(|S| + |T|)^2 / |T| \geq (|S| + |T|)^2 / (|S| + |T|) = |S| + |T|$, whereby}
  &\leq \lim_{n \to \infty}\sum_{\substack{S, T \subseteq V \\ S \cap T = \emptyset \\ |T| \geq 2 \text{ or } |S| + |T| < n}}\exp\left(-\frac{1}{96}C_{\reg, 1}^2 \log n \cdot (|S| + |T|)\right) \nonumber \\
  &\leq \lim_{n \to \infty} \sum_{\substack{S, T \subseteq V \\ |S| + |T| > 0}}\exp\left(-\frac{1}{96}C_{\reg, 1}^2 \log n \cdot (|S| + |T|)\right) \nonumber
  \intertext{and introducing scalar variables $a = |S|$ and $b = |T|$ and grouping according to these values, we further bound}
  &= \lim_{n \to \infty} \sum_{\substack{a, b \in \{0, \dots, n\} \\ a + b > 0}}\binom{n}{a}\binom{n}{b}\exp\left(-\frac{1}{96}C_{\reg, 1}^2 \log n \cdot (a + b)\right). \nonumber
  \intertext{To finish, we use that $\binom{n}{a} \leq \exp(a\log n)$ and $\binom{n}{b} \leq \exp(b \log n)$, and there are at most $(n + 1)^2$ terms in the outer sum, whereby}
  &\leq \lim_{n \to \infty} \sum_{\substack{a, b \in \{0, \dots, n\} \\ a + b > 0}}\exp\left(-\frac{1}{96}C_{\reg, 1}^2 \log n \cdot (a + b) + \log n \cdot (a + b)\right) \nonumber \\
  &\leq \lim_{n \to \infty} \exp\left(-\left(\frac{1}{96}C_{\reg, 1}^2 - 1\right) \log n + 2\log(n + 1)\right).
\end{align*}
Setting $C_{\reg, 1}$ sufficiently large will make the remaining limit equal zero, giving the result.

\sbsctn{Recovery guarantees for the SDP relaxation}

\label{sdpTim}
In the following we state a set of sufficient conditions under which a well-known SDP relaxation of the min-bisection problem
recovers the planted bisection under the SBM. This result follows immediately from Lemma~3.13 of~\cite{Ban18}
(see~\cite{Ban18} for further details about the SDP relaxation). To prove this result we make use of Bernstein's inequality (see Theorem 2.8.4 in~\cite{VerBookHDP}).


\begin{proposition}\label{sdp}
    Suppose $p = p(n) \in (0,1)$, $q = q(n) \in (0, 1)$ are such that the following conditions hold:
    \begin{enumerate}
    \item $\frac{\log n}{3n} < p < \frac{1}{2}$ for sufficiently large $n$.
    \item There exists $\epsilon > 0$ such that $p - q \geq (12 + \epsilon)\sqrt{p \cdot \frac{\log n}{n}}$ for sufficiently large $n$.
    \end{enumerate}
    Then with high probability, the SDP relaxation recovers the planted bisection in a graph drawn from $\mathcal{G}_{n, p, q}$.
\end{proposition}

\begin{prf}
    Define
    \begin{align*}
      \deg_{\tin}(v) &= \begin{cases} N(\{v\}, V_1) & \text{if } v \in V_1 \\ N(\{v\}, V_2) & \text{if } v \in V_2 \end{cases}, \qquad
      \deg_{\tout}(v) = \begin{cases} N(\{v\}, V_2) & \text{if } v \in V_1 \\ N(\{v\}, V_1) & \text{if } v \in V_2 \end{cases}.
    \end{align*}
    Then, by Lemma 3.13 of~\cite{Ban18}, it suffices to show that for any constant $\Delta > 0$, with high probability
    \[ \min_v \left\{ \deg_{\tin}(v) - \deg_{\tout}(v) \right\} \geq \frac{\Delta}{\sqrt{\log n}}\avg\left[ \deg_{\tin}(v_0) - \deg_{\tout}(v_0) \right] \]
    for an arbitrary fixed node $v_0$ (the expectation on the right-hand side does not depend on this choice).
    Noting that
    \[ \avg\left[ \deg_{\tin}(v_0) - \deg_{\tout}(v_0) \right] = \frac{p - q}{2}n - p, \]
    it suffices to show the weaker statement that, again for any constant $\Delta > 0$, with high probability
    \[ \min_v \deg_{\tin}(v) - \max_v \deg_{\tout}(v) \geq \Delta \frac{p - q}{2} \frac{n}{\sqrt{\log n}}. \]

    Note that for any $v$, $\deg_{\tin}(v)$ is a sum of $\frac{n}{2} - 1$ random Bernoulli variables with mean $p$, and therefore with variance $p(1 - p) \leq p$.
    Therefore, using Bernstein's inequality,
    \begin{align*}
      \prob\left[\deg_{\tin}(v) - \avg[\deg_{\tin}(v)] \leq -t\sqrt{pn\log n}\right]
      &\leq \exp\left(-\frac{\frac{1}{2}t^2 pn\log n}{\frac{1}{2}pn + \frac{2}{3}\sqrt{pn\log n}}\right) \\
      &\leq \exp\left(-\frac{t^2 \sqrt{pn}\log n}{\sqrt{pn} + \frac{4}{3}\sqrt{\log n}}\right) \\
      &\leq \exp\left(- \frac{1}{5} t^2 \log n \right).
    \end{align*}
    Thus since $\avg [\deg_{\tin}(v)] = p(\frac{n}{2} - 1)$, taking $t = 3$, and using a union bound over $v$, it follows that, with high probability,
    \[ \min_{v} \deg_{\tin}(v) \geq \frac{1}{2}pn - 3\sqrt{pn \log n}. \]
    By a symmetric argument for $d_{\tout}(v)$, with high probability,
    \[ \max_{v} \deg_{\tout}(v) \leq \frac{1}{2}qn + 3\sqrt{pn \log n}. \]
    Letting $t = \max\{t_{\tin}, t_{\tout}\}$, we find that, with high probability,
    \[ \min_v \deg_{\tin}(v) - \max_v \deg_{\tout}(v) \geq \frac{p - q}{2}n - 6\sqrt{pn\log n}. \]

    Thus it suffices to show that under the assumptions in the statement, for any $\Delta > 0$ and sufficiently large $n$,
    \[ \frac{p - q}{2}n - 6\sqrt{pn\log n} \geq \frac{p - q}{2}n \cdot \frac{\Delta}{\sqrt{\log n}}. \]
    This follows since, by Assumption (2), there exists $\delta > 0$ such that $6\sqrt{pn\log n} \leq (1 - \delta)\cdot\frac{1}{2}(p - q)n$ for sufficiently large $n$.
\end{prf}

Proposition~\ref{sdp} in particular implies that, if $p = \alpha n^{-\omega}$ and $q = \beta n^{-\omega}$ for any $\omega \in (0, 1)$, the SDP relaxation recovers the planted bisection with high probability provided that $p > q$ (or equivalently $\alpha > \beta$). When $p = \alpha \frac{\log n}{n}$ and $q = \beta \frac{\log n}{n}$, by Lemma 4.11 of~\cite{Ban18}, recovery is guaranteed with high probability provided that $\sqrt{\alpha} - \sqrt{\beta} > \sqrt{2}$.

\sbsctn{Distance distributions in random graphs}

\label{distances}
Consider a graph graph $G \sim \mathcal{G}_{n, p, q}$.
In this section, we present lower and upper bounds on the average distance $\rho_{\tavg}(G)$
and the diameter $\rho_{\max}(G)$ of $G$ in very dense, dense,
and logarithmic regimes. We use these results to prove our non-recovery conditions
given by Theorems~\ref{noRecoveryDense} and~\ref{noRecoverySparse}.
The distance distributions of ER graphs have been studied extensively in the literature ~\cite{Bollobas-2001-RandomGraphs,CL-2001-DiameterSparseGraphs,CL-2002-RandomGraphAverageDistance,Shimizu-2017-DenseRandomGraphDistances}.
To obtain similar results for the SBM, we use the following tool
to relate the distance distributions of SBM graphs to those of ER graphs

\begin{proposition}
    \label{prop:coupling}
    Let $0 < q \leq p < 1$ and let $n$ be an even number.
    Then, there exists a probability distribution over triples of graphs $(G_1, G_2, G_3)$ on a mutual set $V$ of $n$ nodes such that the following conditions hold:
    \begin{enumerate}
    \item The marginal distributions of $G_1$, $G_2$, and $G_3$ are $\mathcal{G}_{n, q}$, $\mathcal{G}_{n, q, p}$, and $\mathcal{G}_{n, p}$, respectively.
    \item With probability 1, $G_1$ is a subgraph of $G_2$ and $G_2$ is a subgraph of $G_3$.
        Consequently, for all $i, j \in V$, $\rho_{G_1}(i, j) \geq \rho_{G_2}(i, j) \geq \rho_{G_3}(i, j)$.
    \end{enumerate}
\end{proposition}
\begin{prf}
    We describe a procedure for sampling the three graphs together.
    Let $V$ be their mutual node set, and fix a uniformly random bisection $V = V_1 \sqcup V_2$ with $|V_1| = |V_2| = \frac{n}{2}$.

    The key observation is the following.
    Let $\mathrm{Ber}(a)$ denote the Bernoulli distribution with probability $a$.
    Then, if $X \sim \mathrm{Ber}(a)$ and $Y \sim \mathrm{Ber}(b)$ are sampled independently, then the binary OR of $X$ with $Y$ has the law $\mathrm{Ber}(a + b - ab)$.
    Now, define $a = \frac{p - q}{1 - q}$.
    Since $q + \frac{p - q}{1 - q} - q \cdot \frac{p - q}{1 - q} = p$, we have that if $X \sim \mathrm{Ber}(q)$ and $Y \sim \mathrm{Ber}(a)$, then the binary OR of $X$ with $Y$ has law $\mathrm{Ber}(p)$.

    We now describe a procedure for sampling the desired triple of graphs.
    First, sample $G_1 \sim \mathcal{G}_{n, q}$.
    Next, for each $i \in V_1$ and $j \in V_2$, sample $X_{ij} \sim \mathrm{Ber}(a)$, independently of $G_1$.
    Let $G_2$ have all of the edges of $G_1$, and also an edge between $i$ and $j$ if $X_{ij} = 1$.
    By the previous observation, the law of $G_2$ is $\mathcal{G}_{n, q, p}$.

    Similarly, for each $i, j \in V_1$ and $i, j \in V_2$, sample $Y_{ij} \sim \mathrm{Ber}(a)$, independently of $G_1$ and $G_2$.
    Let $G_3$ have all of the edges of $G_2$, and also an edge between $i$ and $j$ if $Y_{ij} = 1$.
    By the previous observation, the law of $G_3$ is $\mathcal{G}_{n, p}$.
\end{prf}

\smallskip

We now briefly review the existing results on the distance distributions of ER graphs.
It will be convenient for us to give these results in slightly different form than in the original references, so we derive the statements we will need from prior work below.
Below, for $x > 0$, we write $\lfloor x \rfloor$ for the ``floor function'' or greatest integer not exceeding $x$, and $\{x\} \colonequals x - \lfloor x \rfloor$ for the ``fractional part'' of $x$.

\begin{proposition}\label{prop:diameter-er-dense}
    Suppose $p = \alpha n^{-\omega}$ for some $\omega \in [0, 1)$ and $\alpha > 0$, with $\alpha < 1$ if $\omega = 0$.
    Let $G \sim \mathcal{G}_{n, p}$.
    Then,
    \[ \lim_{n \to \infty}\prob\left[\rho_{\max}(G) = 1 + \left\lfloor \frac{1}{1 - \omega}\right\rfloor\right] = 1. \]
\end{proposition}
\begin{prf}
    First, suppose $\omega \in [0, \frac{1}{2})$.
    Then, our task is to show that $\rho_{\max} = 2$ with high probability.
    By Corollary 10.11(i) of~\cite{Bollobas-2001-RandomGraphs}, it suffices to check that $p^2 n - 2\log n \to \infty$ and $n^2(1 - p) \to \infty$.
    The latter clearly holds for any $\omega > 0$, and holds for $\omega = 0$ since in that case $p = \alpha \in (0, 1)$.
    For the former, we have $p^2n = \alpha n^{1 - 2\omega}$ with $\alpha > 0$ and $\omega < \frac{1}{2}$, and the result follows.

    Next, suppose $\omega \in [\frac{1}{2}, 1)$.
    Define $d \colonequals 1 + \lfloor \frac{1}{1 - \omega}\rfloor \geq 3$.
    By Corollary 10.12(i) of \cite{Bollobas-2001-RandomGraphs}, it suffices to check that $d^{-1}\log n - 3\log\log n \to \infty$, $p^dn^{d - 1} - 2\log n \to \infty$, and $p^{d - 1}n^{d - 2} - 2\log n \to -\infty$.
    The first condition follows since $d$ is a constant.
    For the second condition, we calculate
    \[ p^dn^{d - 1} = \alpha^d n^{-\omega(1 + \lfloor \frac{1}{1 - \omega}\rfloor) + \lfloor \frac{1}{1 - \omega}\rfloor}. \]
    Manipulating the exponent, we have
    \begin{align*}
      -\omega(1 + \lfloor \frac{1}{1 - \omega}\rfloor) + \lfloor \frac{1}{1 - \omega}\rfloor
      &= -\omega + (1 - \omega)\lfloor \frac{1}{1 - \omega}\rfloor \\
      &= (1 - \omega)\left(1 - \left\{\frac{1}{1 - \omega}\right\}\right).
    \end{align*}
    We have $(1 - \omega) \in (0, 1)$ and $1 - \{\frac{1}{1 - \omega}\} \in (0, 1]$, so $p^dn^{d - 1} = \alpha^d n^{\delta}$ for some $\delta > 0$, so the second condition holds.
    Finally, for the third condition we have
    \[ p^{d - 1}n^{d - 2} = \alpha^{d - 1}n^{-\omega(1 + \lfloor \frac{1}{1 - \omega}\rfloor) + \lfloor \frac{1}{1 - \omega}\rfloor - (1 - \omega)} = \alpha^{d - 1}n^{-(1 - \omega)\{\frac{1}{1 - \omega}\}} = \alpha^{d - 1}n^{-\delta^{\prime}} \]
    for some $\delta^{\prime} \geq 0$.
    Thus the third condition also holds, and the result follows.
\end{prf}

\begin{proposition}
    Suppose $p = \alpha \frac{\log n}{n}$ for some $\alpha > 1$.
    Let $G \sim \mathcal{G}_{n, p}$.
    Then, for any $\epsilon > 0$,
    \[ \lim_{n \to \infty}\prob\left[\rho_{\max}(G) \in \left((1 - \epsilon)\frac{\log n}{\log \log n}, (1 + \epsilon)\frac{\log n}{\log \log n}\right)\right] = 1. \]
\end{proposition}
\begin{prf}
    We note that the asymptotic $\frac{\log n}{\log pn} \sim \frac{\log n}{\log \log n}$ holds as $n \to \infty$, since $\log pn = \log \alpha + \log\log n$.
    Thus it suffices to show that
    \[ \lim_{n \to \infty}\prob\left[\rho_{\max}(G) \in \left((1 - \epsilon)\frac{\log n}{\log pn}, (1 + \epsilon)\frac{\log n}{\log pn}\right)\right] = 1. \]
    We explicitly make this trivial rewriting to put our result in the form usually used in the literature.
    The result in this form then follows directly from Theorem 4 of~\cite{CL-2001-DiameterSparseGraphs}.
\end{prf}

\begin{proposition}[Theorem 1 of \cite{Shimizu-2017-DenseRandomGraphDistances}]
    \label{prop:moderately-dense-avg-distance}
    Suppose $p = \alpha n^{-\omega}$ for some $\omega \in (0, 1)$ and $\alpha > 0$.
    Define the quantity
    \[ \mu = \mu(\omega, \alpha) = \left\lceil\frac{1}{1 - \omega}\right\rceil + \text{\emph{\textbf{1}}}\left\{\frac{1}{1 - \omega} \in \mathbb{N}\right\}\exp\left(-\alpha^{\frac{1}{1 - \omega}}\right). \]
    Let $G \sim \mathcal{G}_{n, p}$. Then, for any $\epsilon > 0$,
    \[ \lim_{n \to \infty}\prob\left[\rho_{\tavg}(G) \in \left((1 - \epsilon)\mu, (1 + \epsilon)\mu\right)\right] = 1. \]
\end{proposition}

\begin{proposition}[Theorem 1 of \cite{CL-2002-RandomGraphAverageDistance}]
    \label{prop:sparse-avg-distance}
    Suppose $p = \alpha\frac{\log n}{n}$ for some $\alpha > 1$.
    Let $G \sim \mathcal{G}_{n, p}$.
    Then, for any $\epsilon > 0$,
    \[ \lim_{n \to \infty}\prob\left[\rho_{\tavg}(G) \in \left((1 - \epsilon)\frac{\log n}{\log pn}, (1 + \epsilon)\frac{\log n}{\log pn}\right)\right] = 1. \]
\end{proposition}

\smallskip

Utilizing the above results on distance distribution of ER graphs together with Proposition~\ref{prop:coupling},
in the following we obtain similar results for the SBM.

\begin{proposition}\label{prop:dist-dn-sbm-very-dense}
    Suppose $0 < q < p < 1$ are constants.
    Let $G \sim \mathcal{G}_{n, p, q}$.
    Then, for any $\epsilon > 0$,
    \[ \lim_{n \to \infty} \prob\left[(1 - \epsilon)\left(2 - \frac{p + q}{2}\right) \leq \rho_{\tavg}(G) \leq \rho_{\max}(G) = 2\right] = 1. \]
\end{proposition}
\begin{prf}
    By Proposition~\ref{prop:coupling} and Proposition~\ref{prop:diameter-er-dense}, with high probability $\rho_{\max}(G) = 2$ when $G \sim \mathcal{G}_{n, p, q}$.
    On this event, every adjacent pair of nodes in $G$ has distance 1, and every non-adjacent pair of nodes has distance 2.
    Therefore, we may explicitly calculate the average distance: when $\rho_{\max}(G) = 2$, then
    \[ \rho_{\tavg}(G) = \frac{2}{n^2}\left(1 \cdot |E| + 2 \cdot \left(\binom{n}{2} - |E|\right)\right) = 2\left(1 - \frac{1}{n}\right) - \frac{2}{n^2}|E|. \]
    By Hoeffding's inequality, for any fixed $\delta > 0$, with high probability $|E| \leq (1 + \delta) \cdot \frac{p + q}{2}\binom{n}{2}$ (where the second factor on the right-hand side is $\avg |E|$).
    The result then follows after substituting and choosing $\delta$ sufficiently small depending on the given $\epsilon$.
\end{prf}

\begin{proposition}
    \label{prop:dist-dn-sbm-dense}
    Suppose $p = \alpha n^{-\omega}$ and $q = \beta n^{-\omega}$ for $\alpha, \beta > 0$ and $\omega \in (0, 1)$ such that $\frac{1}{1 - \omega} \notin \mathbb{N}$.
    Let $G \sim \mathcal{G}_{n, p, q}$.
    Then, for any $\epsilon > 0$,
    \[ \lim_{n \to \infty} \prob\left[(1 - \epsilon)\left\lceil\frac{1}{1 - \omega}\right\rceil \leq \rho_{\tavg}(G) \leq \rho_{\max}(G) = \left\lceil\frac{1}{1 - \omega}\right\rceil\right] = 1. \]
\end{proposition}
\begin{prf}
The result follows from using Proposition~\ref{prop:coupling} to compare to ER graphs, and then applying Proposition~\ref{prop:moderately-dense-avg-distance} and Proposition~\ref{prop:diameter-er-dense}.
\end{prf}

\begin{proposition}
    \label{prop:dist-dn-sbm-exc-dense}
    Suppose $p = \alpha n^{-\omega}$ and $q = \beta n^{-\omega}$ for $\alpha > \beta > 0$ and $\omega \in (0, 1)$ such that $\frac{1}{1 - \omega} \in \mathbb{N}$.
    Let $G \sim \mathcal{G}_{n, p, q}$.
    Then, for any $\epsilon > 0$,
    \[ \lim_{n \to \infty} \prob\left[(1 - \epsilon)\left(\frac{1}{1 - \omega} + \exp\left(-\alpha^{\frac{1}{1 - \omega}}\right)\right) \leq \rho_{\tavg}(G) \leq \rho_{\max}(G) = \frac{1}{1 - \omega} + 1\right] = 1. \]
\end{proposition}
\begin{prf}
The result follows from using Proposition~\ref{prop:coupling} to compare to ER graphs, and then applying Proposition~\ref{prop:moderately-dense-avg-distance} and Proposition~\ref{prop:diameter-er-dense}.
\end{prf}

\begin{proposition}
    \label{prop:dist-dn-sbm-sparse}
    Suppose $p = \alpha \frac{\log n}{n}$ and $q = \beta\frac{\log n}{n}$ for $\alpha > \beta > 0$ satisfying $\frac{\alpha + \beta}{2} > 1$.
    Let $G \sim \mathcal{G}_{n, p, q}$.
    Then, for any $\epsilon > 0$,
    \[ \lim_{n \to \infty}\prob\left[(1 - \epsilon)\frac{\log n}{\log\log n} \leq
            \rho_{\tavg}(G) \leq \rho_{\max}(G) \leq (1 + \epsilon)\frac{\log n}{\log\log n}\right] = 1. \]
\end{proposition}
\begin{prf}
    The middle inequality holds for any graph $G$.
    For the left inequality, we use Proposition~\ref{prop:coupling} to construct a graph $G^{\prime} \sim \mathcal{G}_{n, p}$ coupled to $G$ such that $G$ is a subgraph of $G^{\prime}$.
    Then, $\rho_{\tavg}(G) \geq \rho_{\tavg}(G^{\prime})$, and the result then follows by Proposition~\ref{prop:sparse-avg-distance}.
    The right inequality follows from a standard result on the diameter of the SBM, which may be obtained by repeating the analysis of the size of the neighborhood of a node given in Section 4.5.1 of \cite{Abbe18} in the case of logarithmic average degree rather than constant average degree.
\end{prf}

\paragraph{Acknowledgements} The authors would like to thank Afonso Bandeira for fruitful discussions about community detection and the recovery properties of SDP relaxations.


\bibliographystyle{plain}

\bibliography{biblio}

\end{document}